\documentclass[10pt,twoside]{article}

\makeatletter
\evensidemargin \oddsidemargin
\makeatother

\usepackage{latexsym,amssymb,graphicx,color,mathrsfs}

\usepackage{amsmath}

\usepackage{mathrsfs} 
\usepackage{multicol}
\usepackage[shortlabels]{enumitem}

\newcommand{\comment}[1]{}
\newcommand{\R}{{\mathbb R}}  

\newtheorem{theorem}{Theorem}
\newtheorem{itlemma}{Lemma}[section] 
\newtheorem{itproposition}[itlemma]{Proposition}
\newtheorem{itcorollary}[itlemma]{Corollary}
\newtheorem{itremark}[itlemma]{Remark}
\newtheorem{itdefinition}[itlemma]{Definition}
\newtheorem{itexample}[itlemma]{Example}

\newenvironment{lemma}{\begin{itlemma}\rm}{\end{itlemma}} 
\newenvironment{remark}{\begin{itremark}\rm}{\end{itremark}} 
\newenvironment{corollary}{\begin{itcorollary}\rm}{\end{itcorollary}}
\newenvironment{proposition}{\begin{itproposition}\rm}{\end{itproposition}}
\newenvironment{definition}{\begin{itdefinition}\rm}{\end{itdefinition}}
\newenvironment{example}{\begin{itexample}\rm}{\end{itexample}}

\newcommand{\be}[1]{\begin{equation}\label{#1}}
\newcommand{\ee}{\end{equation}}
\newcommand{\bl}[1]{\begin{lemma}\label{#1}}
\newcommand{\ble}[1]{\begin{lemmaex}\label{#1}}
\newcommand{\br}[1]{\begin{remark}\label{#1}}
\newcommand{\bt}[1]{\begin{theorem}\label{#1}}
\newcommand{\bd}[1]{\begin{definition}\label{#1}}
\newcommand{\bp}[1]{\begin{proposition}\label{#1}}
\newcommand{\bc}[1]{\begin{corollary}\label{#1}}
\newcommand{\bfact}[1]{\begin{fact}\label{#1}}
\newcommand{\ber}[1]{\begin{exercise}\label{#1}}
\newcommand{\bex}[1]{\begin{example}\label{#1}}
\newcommand{\bem}[1]{\begin{example}\label{#1}}  
\newcommand{\ec}{\mybox\end{corollary}}
\newcommand{\efact}{\mybox\end{fact}}
\newcommand{\eer}{\mybox\end{exercise}}
\newcommand{\eex}{\mybox\end{example}}
\newcommand{\eem}{\mybox\end{example}}
\newcommand{\el}{\mybox\end{lemma}}
\newcommand{\ele}{\mybox\end{lemmaex}}
\newcommand{\er}{\mybox\end{remark}}
\newcommand{\et}{\qed\end{theorem}}
\newcommand{\ed}{\mybox\end{definition}}
\newcommand{\ep}{\mybox\end{proposition}}
\newcommand{\epr}{\end{proof}}
\newcommand{\bpr}{\begin{proof}}

\newcommand{\ecs}{\end{corollary}}
\newcommand{\eers}{\end{exercise}}
\newcommand{\eexs}{\end{example}}
\newcommand{\eems}{\end{example}}
\newcommand{\els}{\end{lemma}}
\newcommand{\eles}{\end{lemmaex}}
\newcommand{\ers}{\end{remark}}
\newcommand{\ets}{\end{theorem}}
\newcommand{\eds}{\end{definition}}
\newcommand{\eps}{\end{proposition}}
\newcommand{\halmos}{\rule{1ex}{1.4ex}}
\newcommand{\mybox}{\hfill $\Box$} 
\newcommand{\beq}{\begin{eqnarray}}
\newcommand{\eeq}{\end{eqnarray}}
\newcommand{\beqn}{\begin{eqnarray*}}
\newcommand{\eeqn}{\end{eqnarray*}}
\newcommand{\bi}{\begin{itemize}}
\newcommand{\ei}{\end{itemize}}
\newcommand{\ben}{\begin{enumerate}}
\newcommand{\een}{\end{enumerate}}

\newenvironment{proof}{\noindent {\em Proof}.\ }{\hspace*{\fill}$\halmos$\medskip}

\newcommand{\qed}{\hfill \halmos} 

\newcommand{\kk}{{\mathcal K}}
\newcommand{\nn}{{\mathcal N}}
\newcommand{\kl}{{\cal KL}}
\newcommand{\ki}{{\cal K}_\infty}
\newcommand{\norm}[1]{\left\Vert#1\right\Vert}
\newcommand{\normo}[1]{{\left\Vert#1\right\Vert}_{\xo}}

\newcommand{\abs}[1]{{\left\vert#1\right\vert}}
\newcommand{\normx}[1]{{\left\Vert#1\right\Vert}_{\xo}}
\newcommand{\enorm}[1]{\left| #1\right|}

\newcommand{\rref}[1]{(\ref{#1})}
\newcommand{\edo}{\end{document}}

\newcommand{\iss}{{\sc iss}}

\newcommand{\ios}{{\sc ios}}
\newcommand{\siios}{{\sc si}-{\ios}}
\newcommand{\olios}{{\sc ol}-{\ios}}
\newcommand{\ubibs}{{\sc ubibs}}
\newcommand{\ros}{{\sc ros}}

\newcommand{\rfc}{{\sc rfc}}
\newcommand{\ol}{{\sc ol}}
\newcommand{\olgs}{{\sc ol-gs}}
\newcommand{\gaos}{{\sc gaos}}
\newcommand{\rgaos}{{\sc rgaos}}
\newcommand{\olrgaos}{{\sc ol-rgaos}}
\newcommand{\sirgaos}{{\sc si-rgaos}}
\newcommand{\oag}{{\sc o-ag}}
\newcommand{\ogs}{{\sc o-gs}}
\newcommand{\gs}{{\sc gs}}

\newcommand{\xco}{{\cal X}}
\newcommand{\xc}{{\cal \xco}^n}

\newcommand{\ve}{{\varepsilon}}

\newcommand{\vf}{\varphi}

\newcommand{\xo}{{\cal X}}

\newcommand{\Uu}{{\cal U}}
\newcommand{\Aa}{{\cal A}}

\newcommand{\xu}{{x_{t}^{\xi, u}}}

\newcommand{\cN}{{\mathcal N}}

\newcommand{\sigmay}{\sigma_y}

\newcommand{\at}{\overline{\alpha}}
\newcommand{\ab}{\underline{\alpha}}

\newcommand{\bn}{\left\|}

\newcommand{\enx}{\right\|_\mathcal{X}}

\newcommand{\MO}{{\mathcal M_{\Omega}}}  
\newcommand{\Nn}{{\mathcal N}}
\newcommand{\MmU}{{\mathcal M_{\Uu}}}
\newcommand{\Mm}{{\mathcal M}}
\newcommand{\Rr}{{\mathcal R}}
\newcommand{\Zz}{{\mathcal Z}}
\newcommand{\rcomment}{\textcolor{red}}

\begin{document}
\thispagestyle{empty}
\setcounter{page}{1}

\noindent
{\footnotesize {\rm Submitted to\\[-1.00mm]
{\em Dynamics of Continuous, Discrete and Impulsive Systems}}\\[-1.00mm]
http:monotone.uwaterloo.ca/$\sim$journal} $~$ \\ [.3in]

\begin{center}
{\large\bf Remarks on Different Notions on\\
Output Stability for Nonlinear Delay Systems}

\vskip.20in

Hasala Gallolu Kankanamalage$^{1}$,\ 
Yuandan Lin$^{2}$,\ and\ Yuan Wang$^{2}$ \\[2mm]
{\footnotesize
$^{1}$Department of Mathematics, Roger Williams University,\\ Bristol, RI 02809, USA.\\[5pt]
$^{2}$Department of Mathematical
  Sciences, Florida Atlantic University,\\ Boca Raton, FL 33431, USA\\
Corresponding author email:
{\tt hgallolu@rwu.edu}
}
\end{center}

{\footnotesize
\noindent
{\bf Abstract.}
Motivated by the regulator theory and adaptive controls, several notions on output stability in the framework of input-to-state stability (\iss) were introduced for finite-dimensional systems. It turned out that these output stability notions are intrinsically different, reflecting different manners of how state variables may affect the transient behavior of output variables.
In this work, we consider these output stability properties for delay systems. Our main objective is to illustrate how the various notions are related for delay systems and to provide the Razumikhin criteria for the output stability properties. The main results are also critical in developing the converse Lyapunov theorems of the output stability properties for delay systems.
\\[3pt]
}


\section{Introduction}

Since introduced in \cite{sontag-coprime},
the \iss\ theory has become a natural framework for robust stability analysis for systems
affected by disturbances and uncertainties.
Motivated by applications such as tracking and regulation, adaptive control, and observer design, a set of notions on robust output stability were introduced in
\cite{sw-ios-scl} and \cite{sw-ios-siam}. These notions become
especially relevant in a situation when it is only possible or necessary to investigate
stability properties for a partial set of state variables, or more generally, functions of
state variables.
Such a situation frequently
arises from applications of regulation, adaptive control design, etc.
Stability properties such as stability
respect to sets (e.g., closed orbits) and incremental stability can also be phrased as
appropriate output stability properties. Some related properties were
treated in the past literature
as partial stability (see \cite{Vorotnikov-05}), or stability ``in two
measures'' (see \cite{LALI93}, \cite{TP-ESAIM00}, and \cite{Kellet-Dower-CDC12}.)
Other developments and applications of output stability properties
can be found in, for instance, \cite{Che-Haddad-2002}, \cite{Fukui-Satoh-cdc20}, \cite{JH-book}, \cite{Karaf-05},  \cite{LWC}, \cite{JoseLuis-Garcia},
\cite{Zimenko-Efimov-Polyaov-Kremlev},  
and notably for delay systems, in \cite{klw-ifac17}, \cite{kar-jiang-book}, \cite{Kar07c}, \cite{Kar07d},  and \cite{POLILU06}.

Roughly speaking, for a system to be robustly output stable, it should be required that,
for any initial state, the magnitude of the output $y(\cdot)$ reduces to a magnitude proportional to the size of the disturbance $u$.
The work \cite{sw-ios-scl} formulated a set of notions associated with the different converging behaviors. These notions may seem to be similar, but they reflect precisely the different converging behavior of the output variables, resulting in substantially different consequences. For instance, the extra output-Lagrange global stability (\olgs) condition imposed on the input-to-output stability (see Section \ref{s-3} for precise definitions) plays a critical role in the converse Lyapunov results, at least for the delay-free systems (see \cite{sw-ios-siam}.)  The \olgs\ property can also lead to progress along the lines of Lyapunov-Krasovskii functionals with point-wise dissipation rates, an important problem specific for delay-systems, see \cite{CHPEMACH17} and \cite{CHGOPE21}. 

Our first contribution is to provide the Razumikhin criteria for the output stability properties. The Razumikhin method is commonly used for stability analysis for delay systems, see for instance,  \cite{Razumikhin56}, \cite{Hale-book}, and \cite{teel-delay}. From the point of view of small-gain theorems, an advantage of the Razumikhin method is that it can convert the task of stability analysis for delay systems to a task of robust stability analysis for delay-free systems by treating the delayed variables as system uncertainties. In the work  \cite{Kar07d}, a Lyapunov-based Razumikhin result was provided for a few types of input-to-output stability properties for time-varying systems with delays.
In this work, we provide some trajectory-based Razumikhin criteria for the different types of output stability properties studied in this paper.The proofs of the Razumikhin results were motivated by the work on small-gain theorems in \cite{TEELtac96}.
Our second contribution is to extend the work and results in \cite{sw-ios-scl}
to time-delay systems. 
There are substantial obstacles when carrying out the analysis for delay systems. For instance,  one of the main results in \cite{eds-david-fc} asserts that a delay-free system is robustly forward complete which means that the reachable set from a bounded set of initial states under bounded control stays bounded in finite time (see also \cite{iss-new}),  yet it remains an open question if this boundedness property still holds for delay systems (see more discussions in \cite{karafyllis2022global}.)
In fact, many of the open questions can be traced to the question whether forward completeness is equivalent to the so called robust forward completeness property, used quite widely in the literature of infinite dimensional or delay systems, see \cite{kar-jiang-book} and \cite{mironchenko2017characterizations}. Following the work
\cite{kar-jiang-book} and \cite{mironchenko2017characterizations}, we adopt the robust forward completion property as one of our basic assumptions.  In the Conclusion section, we also indicate the limitations of our results when comparing to the ones for delay-free systems.

The current work is also a preparation
for developing converse Lyapunov results of the \ios-type
properties for delay systems (see \cite{klw-ifac17}).
Some of these results were reported in a preliminary version \cite{klw-icca17} with some roughly sketched proofs. Since then the work has been developed further. One of the main assumptions of bounded-input-bounded-state in \cite{klw-icca17} is now relaxed to the robust forward completeness condition.  We also treat the stability notions more thoroughly to allow both output maps with delayed state variables and delay-free output maps. 
%
%

The rest of the paper is organized as follows. In Section \ref{s-2}, we go through a preliminary setup of the systems to be considered in this work, together with some brief discussions on some immediate consequence of the robust forward completeness. In Section \ref{s-3} we present definitions on the output stability properties. It is worth noting that with the definitions given for general output maps on possibly delayed state variables, it makes sense to consider the stability notions in terms of ``history norms''.  In Section \ref{s-raz}, we provide the Razumikhin results for each of the output stability properties.
In Section \ref{s-4}, we show how the \ios\ property can be converted to \olios\ under output redefinition.  This will be a critical step for our future work in developing a converse Lyapunov theorem for the \ios\ property (see \cite{sw-ios-siam} for the delay-free case.) In Section \ref{s-5}, we show that if a system has the \ios\ (or \olios) property, then the system admits $\ki$-margin under output feedback to preserve the corresponding output stability when $u\equiv 0$.  The proofs of the main results in Sections \ref{s-raz}, \ref{s-4}, and \ref{s-5} are given in Section \ref{s-6}. We conclude the paper by Section \ref{s-7} with some discussions on some open problems.  For the convenience of references, we give detailed proofs of some preliminary results in the Appendix.

{\it Notations}. Throughout this work, we use $|\cdot|$ to denote the
Euclidean norm in $\R^n$,
and $\|\cdot\|_{I}$ for the $L_\infty$
norm of measurable functions on the interval $I$.
For a subset $\Uu$ of $\R^m$, we use $\Mm_\Uu$ to denote the set of measurable and locally essentially bounded functions from $\R_{\ge 0}$ to $\Uu$, and simply use $\Mm$ when $\Uu=\R^m$.

A continuous function $\sigma: \R_{\geq 0} \to \R_{\geq 0}$ is of {\em class $\Nn$} if $\sigma(0)=0$ and is nondecreasing;
is of {\em class $\kk$} if it is of class $\Nn$ and is strictly increasing; and is of {\em
class $\ki$} if it is of class $\kk$ and onto.  A function $\beta : \R_{\geq 0}
\times \R_{\geq 0} \to \R_{\geq 0}$ is said to be  of {\em class $\kl$}
if for each fixed $t \geq 0$, $\beta(\cdot, t)$ is of class $\kk$, and
for each fixed $s \geq 0$, $\beta(s,t)$ decreases to 0 as $t \to
\infty$.

In this work we consider systems with a finite delay interval $[-\theta, 0]$, $\theta>0$.  Denote $C([-\theta, 0], \R^k)$ by $\xo^k$,
the Banach space of continuous functions from $[-\theta, 0]$ to $\R^k$, equipped with the $C^0$ norm $\normo{\cdot}$. 
For a continuous function $v$ defined on $[-\theta, b]$, define, for each $t\in [0, b]$, the history segment of $v$ at $t$ by
$v_t(s):=v(t+s)$ for $s\in [-\theta, 0]$.

A continuous function $F$ from a normed vector space ${\cal A}$ to a normed vector space ${\cal B}$
is said to be completely continuous if $F$ maps a bounded set to a bounded set, that is,
for each bounded set $S\subseteq {\cal A}$, $F(S)$ is a bounded subset of $\mathcal B$.

A function $F$ from a normed vector space ${\cal A}$ to a normed vector space ${\cal B}$
is said to be Lipschitz on compact sets if for every compact subset $K$ of $\mathcal A$, there is some $M\ge 0$ such that
\be{e-comp}
\norm{F(p)-F(q)}_{\cal B}\le M\norm{p-q}_{\cal A},
\ee
for all $p, q\in K$; $F$ is said to be locally Lipschitz if for every $p_0\in\mathcal A$, there is a neighborhood  $O$ of $p_0$ and some $M\ge 0$ such that \rref{e-comp} holds for all $p, q\in O$;
and $F$ is said to be Lispschitz on bounded sets if for each bounded subset $B$ of $\mathcal A$, there is some $M\ge 0$ such that \rref{e-comp} holds for all $p, q\in B$. The following can be seen readily:
\[
\mbox{Lipschitz on bounded sets $\Rightarrow$ locally Lipschitz $\Rightarrow$ Lipschitz on compact sets}.
\]
When $\mathcal A$ is finite dimensional, all three Lipschitz properties are equivalent.  The following was shown
in \cite{Xu-Lip-2020} for the general case when $\mathcal A$ is not necessarily finite dimensional:
\bl{l-Xu-Lip}
Let $\mathcal A$ and $\mathcal B$ be normed vector spaces. A function $F:\mathcal A\rightarrow\mathcal B$ is locally Lipschitz if and only if $F$ is Lipschitz on compact sets.
\el
It was also shown in \cite{Xu-Lip-2020} by a counterexample that a locally Lipschitz map $F$ is not necessarily Lipschitz on bounded sets.

\section{Preliminaries}\label{s-2}

In this work we consider systems with a finite delay interval $[-\theta, 0]$ as follows:
\begin{subequations}\label{e-syso}
\begin{align}
\dot x(t) &= f(x_t, u(t)), \label{e-sysox}\\
y(t) &= h(x_t), \label{e-sysoo}
\end{align}
\end{subequations}
where for each $t$, $x(t)\in\R^n, \ u(t)\in\R^m$, and $y(t)\in\R^p$. It can be
interesting to consider more general output maps, for instance,  when
$y(t)$ takes values in a normed space ${\cal Y}$ with a norm $\norm{\cdot}_{\cal Y}$.  But in the context of output stability, one can always consider $\norm{y(t)}_{\cal Y}$ as the output.

The map $f: \xo^n\times\R^m \rightarrow \R^n$ is assumed to be completely continuous and locally Lipschitz; and the output map $h:\xc\rightarrow\R^p$ is assumed to be continuous, satisfying the boundedness condition for some $\pi\in \cN$:
\be{bound-h}
\abs{h(\xi)} \leq \pi \left( \bn \xi \enx \right),
\ee
and in particular, $h(0) = 0$. The output map $h$ is said to be delay-free if there exists a continuous map $h_0: \xc\rightarrow\R^p$ such that $h(\xi) = h_0(\xi(0))$.

The inputs $u(\cdot)$ are measurable and locally essentially bounded
functions defined on $\R_{\ge 0}$.  With the assumptions imposed on $f$, the
usual regularity conditions such as existence, uniqueness, and maximum
continuation of solutions hold for system~\rref{e-sysox}. More specifically, for each initial
state $\xi \in\xo^n$ and input function $u(\cdot)\in\Mm$, there is a unique solution of
\rref{e-sysox} defined on its maximum interval $[-\theta, T_{\xi, u})$
which is locally absolutely continuous on $[0, T_{\xi, u})$ and satisfies the equation \rref{e-sysox} a.e. on $[0, T_{\xi, u})$, fulfilling the initial condition $x_0=\xi$ (see \cite[Section 2.6]{Hale-book} and \cite[Section 3--2.2]{kolmanovskii-book}.)
We use
$x(\cdot, \xi, u)$ to denote such a  maximal solution of
\rref{e-sysox}, and accordingly, we use
$x_t^{\xi, u}: [-\theta, 0]\rightarrow\R^n$ to denote
the map induced by
$x(t, \xi, u)$, that is, $x_t^{\xi, u}(s) = x(t+s, \xi, u)$.
We use $y(t, \xi, u)$ to denote the corresponding
output function $h(x_{t}^{\xi, u})$, i.e., $y(t, \xi, u) = h(x_t^{\xi, u})$.  If the output map $h$ is delay-free, then the history segment of $y$, given by $y_t^{\xi, u}\in \xo^p$, is defined for $0\le t < T_{\xi, u}$.

We say that system \rref{e-sysox} is {\it forward complete} if for each $\xi$ and
each $u(\cdot)$, the solution $x(\cdot, \xi, u)$ is defined for all
$t\ge 0$.

One of the main results in \cite{lsw-siam} asserts that for a delay-free system,
if it is forward complete, then trajectories starting from a compact set with bounded inputs stay
in a compact set in finite time (Proposition 5.1 in \cite{lsw-siam}).  It is not clear if this is still the case for delay systems.
In the work \cite[Definition 2.1]{kar-jiang-book} and \cite[Definition 4]{mironchenko2017characterizations}, the following \rfc\ property (also referred as bounded reachability sets in finite time) was defined:
\bd{d-rfc}
A forward complete system \rref{e-sysox} is said to be {\it robustly forward complete} (\rfc) if for every $T, R>0$,
\be{e-rfc}
\sup\left\{\|x_t^{\xi, u}\|:\ \normo{\xi}\le R, \ \norm{u}\le R, \ t\in[0, T]\right\} < \infty.
\ee
\eds
It was shown in \cite{eds-david-fc} that for a delay free system, the \rfc\ property is equivalent to the forward completeness property.  For delay systems, it still remains open if the forward completeness property implies the \rfc\ property.
A few preliminary properties related to the \rfc\ properties that will be used in this work are listed below. (Recall that $\MmU$ denotes the set of measurable and locally essentially bounded functions from $\R_{\ge 0}$ to $\Uu$, and $\Mm$ is used for $\MmU$ when $\Uu=\R^m$.)

\bl{l-rfc}
System \rref{e-sysox} is \rfc\ if and only if
there exist $\chi_i\in\cN$ ($i=1, 2, 3$) and $c\ge 0$ such that the following holds along all trajectories of the system:
\be{e-rfc-e}
\abs{x(t, \xi, u)}
\le \chi_1(t) + \chi_2(\normo{\xi}) + \chi_3(\norm{u}) + c,\qquad
\ee
for all $t\ge 0$, all $\xi\in\xc$, and all $u\in\Mm$.
\el

Assume that system \rref{e-sysox} is forward complete. Define
\[
\Rr^{T}_{\Uu}(S) = \{\eta: \ \eta=x_t^{\xi, u},\ 0\le t\le T,  u\in\MmU,\ \xi\in S\}.
\]
For a general subset $\Aa$ of $\xc$, we use $\overline{\Aa}$ to denote the closure of $\Aa$.

\comment{
\bp{p-ufc} A system as \rref{e-syso} satisfies the \rfc\ property if and only if for each bounded $\Uu\subseteq
\R^m$, each bounded $S \subseteq \xc$, and $T>0$, the set
${\Rr_{\Uu}^{T}(S)}$ is bounded.
\eps
}

\bl{l-comp} Assume that system \rref{e-sysox} is \rfc.
Then the set $\overline{\Rr_{\Uu}^{T}(S)}$ is compact if $S\subseteq\xc$ is compact, $\Uu\subseteq\R^m$ is bounded, and $T>0$.
\el

\bc{c-lip} Assume that system \rref{e-sysox} is \rfc.
\begin{enumerate}
\item[(a)]
Assume that $S\subseteq\xc$ is compact, $\Uu\subseteq\R^m$ is bounded, and $T>0$.  Then the map
$\xi\mapsto x(t, \xi, u)$ is Lipschitz on $S$ uniformly in $t\in [0, T]$ and $u\in\MmU$.
\item[(b)]
Suppose further that $f$ is Lipschitz on bounded sets.  Assume that $S\subseteq\xc$ is bounded, $\Uu\subseteq\R^m$ is bounded, and $T>0$. Then the map $\xi\mapsto x(t, \xi, u)$ is Lipschitz on $S$ uniformly in $t\in [0, T]$ and $u\in\MmU$.\mybox
\end{enumerate}
\ecs

The proofs of Lemmas \ref{l-rfc} and \ref{l-comp} and Corollary \ref{c-lip} will be given in Appendix \ref{s-A1}.

\section{Notions on Input-to-Output Stability Properties}\label{s-3}

In this section we extend the notions on output stability given in \cite{sw-ios-scl} to time-delay systems as in \rref{e-syso}. The \ios\ property was also introduced in \cite{Kar07d} for more general types of systems.  In this work we focus on time-invariant delay systems, and explore two other related notions: \olios\ and \siios\ as defined below.  

\bd{d-ios}(\cite{sw-ios-scl}, \cite{Kar07d}) A forward complete system as in \rref{e-syso} is:
\begin{itemize}
\item \textit{input-to-output stable } (\ios) if there exist $\beta \in
  \kl$ and $\gamma \in \Nn$ such that
\be{ios}
\abs{y(t, \xi, u)}\leq \beta \left( \bn \xi \enx , t \right) +
\gamma \left(\norm{u}\right); 
\ee
\item \textit{output Lagrange input-to-output stable } (\olios) if it
  is \ios\ and, additionally, the following \olgs\ property holds: there exists $\sigma \in \Nn$ such
  that
\be{ol}
\hspace{-3pt}
\abs{y(t, \xi, u)}\leq \max\!\left\{ \sigma \left(\abs{h(\xi)}
\right), \sigma \left(\norm{u} \right) \right\};
\ee
\item \textit{state independent input-to-output stable } (\siios) if
  there exist $\beta \in \kl$ and $\gamma \in \Nn$ such that
\be{siios}
\abs{y(t, \xi, u)} \leq \beta \left(\abs{h(\xi)}, \ t \right)
+ \gamma \left(\norm{u}\right); 
\ee
\end{itemize}
where all the estimates listed above are assumed to hold for all $\xi, u$ and
all $t\ge 0$.
\ed

For a delay-free output map $h$, properties \rref{ol} and \rref{siios} can be too strong (which can be seen from the special case when $h(\xi) = h_0 \left( \xi(0)\right)$.) This motivates one to consider the output history map $Y=H(\xi)$ associated with $h_0$, defined by
\be{e-H}
H(\xi) = \max_{-\theta\le s\le 0}\abs{h_0(\xi(s))}.
\ee
Note that if $h_0$ is locally Lipschitz, then $H:\xc\rightarrow\R$ is Lipschitz on bounded sets. Along the system trajectories, $Y(t, \xi, u) = \|y_t^{\xi,u}\|_\xo$.
By definition, the \olgs\ property \rref{ol} in terms of the output $Y$ becomes
\be{e-newol}
\|y_t^{\xi, u}\|_{\xo}\le\max\left\{\sigma\left(H(\xi)\right), \ \sigma(\norm{u})\right\},  \quad \forall\,t\ge 0,
\ \forall\,\xi\in\xc, \ \forall\,u\in\Mm.
\ee
This is equivalent to the following, where $\sigma_1\in\kk$:
\be{e-newol1}
\abs{y(t, \xi, u)}\le \max \left\{\sigma_1\left(H(\xi)\right), \ \sigma_1(\norm{u})\right\},  \quad \forall\,t\ge 0, \ \forall\,\xi\in\xc, \ \forall\,u\in\Mm.
\ee
We say that a forward complete system with a delay-free output map $h_0$ is \olios\ {\it in the history norm\/} if the system is \ios\ and satisfies the \olgs-property \rref{e-newol1}.
Similarly, a forward complete system with a delay-free output map $y=h_0(\xi)$ is said to be \siios\ in the history norm if there exist $\beta\in\kl$ and $\gamma\in\cN$ such that for all $\xi\in\xc$ and all $u\in\Mm$,
\be{e-newsiios}
\abs{y(t, \xi, u)}\le \beta\left(H(\xi), \ t \right)
+ \gamma \left(\norm{u}\right), \quad \forall\, t \geq 0.
\ee
It should be noticed that when the output map $y=h_0(\xi(0))$ is delay-free, the \olios\ and \siios\ properties in the history norm are more natural to be utilized in applications than the \olios\ and \siios\ properties directly in terms of $y$. However, for the general case when $h$ depends on $\xi$ over the delay interval $[-\theta, 0]$, the history segment of $y$ at a moment $t<\theta$ is not well defined. This is the main reason that the \ios-properties in Definition~\ref{d-ios} were not formulated by using the history norm of $y$.

\comment{
It can be seen that if a system \rref{e-syso} is \ios, then the
system is zero-input output stable, that is, the
zero-input system
\beq\label{e-syso-0input}
\dot x(t) = f_0(x_t)\  (:=f(x_t, {0})), \
y(t) = h(x_t),
\eeq
satisfies the following output stable property:
\[
\abs{y(t, \xi, 0)} \leq \beta \left( \bn \xi \enx , t
 \right), 
\]
for all $\xi\in\xc$ and all $t\ge 0$.
}

The \olgs\ property \rref{ol} plays a critical role in establishing the converse Lyapunov results for the several output stability properties. The property can also be significant in further work such as point-wise dissipation rate of Lyapunov functions, and possibly lead to progress along the line of uniform output convergence property, see \cite{Ingalls-ACC01} for the delay-free case.
As in the delay-free systems in \cite{sw-ios-scl}, the notion \olios\
can be described by a more compact representation that combines
the \ios\ property together with output Lagrange property \rref{ol}.
To make this work more self-contained, we provide a detailed proof of the result in Appendix~\ref{s-olios}.

\bp{2nd-form-olios}
A forward complete system \rref{e-syso} is \olios\ if and only
if for some $\beta\in\kl$ and $\rho, \gamma\in\Nn$:
\be{e-olios2}
\abs{y(t, \xi, u)}\le\beta\left(\abs{h(\xi)},\
  \frac{t}{1+\rho(\normx{\xi})}\right) + \gamma(\norm{u})
\ee
for all $t\ge 0$, all $\xi\in\xc$ and all $u\in\Mm$.
\eps

An estimate as in \rref{e-olios2} captures both the \ios\ condition
and the Lagrange property  \rref{ol}. It reflects the essential feature of \olios\
property that the output overshoot is dominated by the
initial output value (instead of by the initial value of the state variables
as in the \ios\ case), but the decay rate of the output function
may be slowed down by large initial state values.  For the \siios\ condition,
both the overshoot and the decay rate are determined by the initial value of the output.

By causality, one may use $\gamma(\norm{u}_{[0, t)})$ to replace
$\gamma(\norm{u})$ in the estimates \rref{ios}, \rref{ol}, \rref{siios}, \rref{e-newsiios}
and \rref{e-olios2}.
\comment{not needed?
It is not hard to see the implications \siios\ $\Rightarrow$ \olios\ $\Rightarrow$ \ios.
In the special case when the output function $h$ is given by $h(\xi) = \normo{\xi}$,
all three output stability
notions become the well known ISS property.  However, in the
general case, even for delay-free systems,
none of the reversed implications holds (see \cite{sw-ios-scl}), and
thus, the three notions define different stability properties.
}
Also observe that it results in equivalent definitions if
requirement \rref{ios} is replaced by
\be{ios-max}
\abs{ y(t, \xi, u)} \leq \max\left\{
\beta \left(\normo{\xi}, \ t \right), \
\gamma \left(\norm{u}\right)\right\}, 
\ee
and similarly for the \siios\ expression \rref{siios}.
\comment{
 by
\be{siios-max}
\abs{y(t, \xi, u)} \leq \max\left\{
\beta \left(\abs{y_0}, \ t \right), \
\gamma \left(\norm{u}\right)\right\}. 
\ee
}

We associate with system \rref{e-syso} the following zero-input
system:
\be{e-syso-0input}
\dot x(t) =  f_0(x_t) := f(x_t, 0),\quad
y(t) = h(x_t).
\ee
Note that if a system \rref{e-syso} is \ios, then
the zero-input system \rref{e-syso-0input} is {\it globally asymptotically output stable} (\gaos):
\[
\abs{y(t, \xi, 0)} \leq \beta \left( \bn \xi \enx , t
 \right), 
\]
for all $\xi\in\xc$ and all $t\ge 0$. Similarly, property
\rref{ol} renders to \rref{e-syso-0input} the following output Lagrange stability property:
\[
\abs{y(t, \xi, 0)} \le\sigma_1(\abs{h(\xi)}),\quad\forall\,t\ge 0.
\]

Observe that any of \ios, \olios, and \siios\ implies the following {\it output asymptotic gain} (\oag) property: there exists some $\gamma\in\Nn$ such that the following holds for all trajectories:
\be{e-oag}
\limsup_{t\rightarrow\infty}\abs{y(t, \xi, u)}\le \gamma(\norm{u}).
\ee
The \ios\ property \rref{ios} implies the following {\it output global stability} (\ogs) property: there exist $\sigma\in\ki$ and $\gamma\in\Nn$, the following holds for all trajectories:
\be{e-ogs}
\abs{y(t, \xi, u)}\le \max\{\sigma(\normo{\xi}), \ \gamma(\norm{u})\},\quad\forall\,t\ge 0.
\ee
The \olios\ (and consequently the \siios) property implies the (\olgs) property as defined by \rref{ol}.  Hence, we get the following:
\bp{p-oag}
Assume system \rref{e-syso} is forward complete.
\begin{itemize}
    \item[(a)] If the system is \ios, then it satisfies the \ogs\ and the \oag\ property.
    \item[(b)] If the system is \olios, then it satisfies the \olgs\ and the \oag\ property.\mybox
\end{itemize}
\eps

The \oag, \ogs,  and \olgs\ properties were used for output stability analysis in the past literature, see for instance, \cite{TEELtac96} and \cite{Ingalls-ACC01} for delay-free systems, and \cite{PTM-2006} and \cite{twj-dcdis} for delay systems. It should be noted that the converse of statement (a) fails even for delay-free systems. It still remains open if the converse of statement (b) holds for delay-system, though it was shown to be true for delay-free system, see \cite{Ingalls-ACC01}.

\section{The Razumikhin Method}\label{s-raz}
In this section we provide the Razumikhin results on \ios, \olios, and \siios\ for systems with continuous delay-free output maps.

\bd{d-gs} Consider a system as in \rref{e-sysox}. We say that the system satisfies:
\begin{enumerate}
    \item the uniformly bounded-input bounded-state (\ubibs) property if there exist some $\sigma\in\ki$, $\mu\in\Nn$, and $c>0$ such that the following holds for all $\xi\in\xc$, all $u\in\Mm$:
\be{e-ubibs}
\abs{x(t, \xi, u)}\le \sigma(\normo{\xi}) + \gamma(\norm{u}) + c,
\quad\forall\,t\ge 0;
\ee
\item the global stability (\gs) property if the \ubibs\ property \rref{e-ubibs} holds with $c=0$.~\mybox
\end{enumerate}
\eds

Consider a system as in \rref{e-sysox} with a continuous delay free output map $y=h_0(x)$ so that along
trajectories of the system, 
$y(t, \xi, u) = h_0(x(t, \xi, u))$.  In \cite{Kar07d}, a Razumikhin theorem was provided for such a system in terms of Lyapunov functions.
Below we provide some trajectory-based Razumikhin results. Recall that when the output map $h_0$ is delay-free, the associated output history map $H:\xc\rightarrow\R_{\ge 0}$ is defined by \rref{e-H} so that along trajectories of the system, $H(x_t^{\xi,u})=\|y_t^{\xi, u}\|_{\xo}$.
For simplicity, we use $y_0$ to denote $y_0^{\xi, u}$ when the initial state is clear in the context.
\bp{p-raz-ios}
Consider system \rref{e-sysox} with a continuous delay-free output map $y=h_0(x)$.  
\begin{enumerate}
\item Suppose that the system \rref{e-sysox} is \gs, and assume that there exist $\beta\in\kl$,  $\kappa, \gamma\in\Nn$ with $\kappa(s)<s$ for all $s>0$ such that for all $\xi\in\xc, \ u\in\Mm$, the following holds along all trajectories:
\be{e-raz-ios}
\abs{y(t)}\le\max\left\{\beta\left(\normo{\xi},\  t\right), \ \kappa\left(\max_{\tau\in [-\theta, t]}\abs{y(\tau)} \right),\ \gamma(\norm{u})\right\},
\ee
for all $t\ge 0$. Then the system is \ios.
\item Suppose that the system \rref{e-sysox} is \ubibs, and assume that there exist $\beta\in\kl$, $\rho, \kappa, \gamma\in\Nn$ with $\kappa(s)< s$ for all $s>0$ such that for all $\xi\in\xc, \ u\in\Mm$, the following holds along all trajectories:
\beq
\abs{y(t)} &\le& \max\left\{\beta\left(\normo{y_0}, \ \frac{t}{1+\rho(\normo{\xi})}\right), \right. \nonumber\\ 
& & \qquad\qquad \left.
\kappa\left(\max_{\tau\in [-\theta, t]}\abs{y(\tau)} \right),\ \gamma(\norm{u})\right\},
\label{e-raz-ol}
\eeq
for all $t\ge 0$. Then the system is \olios\ in the history norm.
\item Suppose that system \rref{e-sysox} is forward complete, and assume
there exist $\beta\in\kl$,  $\kappa, \gamma\in\Nn$ with $\kappa(s)<s$ for all $s>0$ such that  for all $\xi\in\xc, \ u\in\Mm$, the following holds along all trajectories:
\be{e-raz-siios}
\abs{y(t)}\le\max\left\{\beta(\normo{y_0}, \ t), \ \kappa\left(\max_{\tau\in [-\theta, t]}\abs{y(\tau)} \right),\ \gamma(\norm{u})\right\},
\ee
for all $t\ge 0$. Then the system is \siios\ in the history norm.
\end{enumerate}
\eps

In the special case when $y=x(t)$, that is, when $h(\xi)=\xi(0)$, conditions \rref{e-raz-ios}, \rref{e-raz-ol}, and \rref{e-raz-siios} are all equivalent. For this case, the assumption of forward completeness for Statement 3 is not needed, and by assuming that \rref{e-raz-siios} holds on the maximum intervals, one can show that the system is forward complete. See \cite[Proposition 2.3]{twj-dcdis} for more details of this special case.

\bex{ex-raz}
To illustrate the meaning of Proposition \ref{p-raz-ios}, consider the parameterized system
\be{e-ex2}
\dot x(t) = \frac{-x(t-\theta) + u(t)}{1+p^2},
\ee
where $p$ is an unknown parameter that takes a constant value, and $\theta\in [0, 1]$. To consider the stability property of the system, we consider the augmented system, where $p$ is treated as a state variable:
\be{e-ex21}
\begin{array}{ccl}
\dot p(t) &=& 0,\\[2mm]
\dot x(t) &=& \frac{\displaystyle{-x(t-\theta)+u(t)}}{\displaystyle{1+p(t)^2}},
\end{array}
\ee
with $y(t) = x(t)$ as output, i.e., $h(\xi) = \xi(0)$, a delay-free output. Solutions of the $p$-subsystem are all constants.  Consider the function
\[
V(p, x) = \frac{x^2}{2}.
\]
Along the trajectories of \rref{e-ex21}, it holds that
\[
\frac{d}{dt}V(p(t), x(t))
=\frac{-x(t)x(t-\theta) + x(t)u(t)}{1+p^2(t)}=\frac{-x(t)x(t-\theta) + x(t)u(t)}{1+p^2_0},
\]
where $p_0=p(0)$.
Note that for any $t>0$, there exists some $\hat t\in (t-\theta, t)$ such that
\beqn
x(t-\theta) &=& x(t) - (x(t) - x(t-\theta))
= x(t) - \dot x(\hat t)\theta\\
&=&x(t) - \frac{(-x(\hat t-\theta) + u(\hat t))\theta}{1+p_0^2}.
\eeqn
Therefore,
\beqn
& &\frac{d}{dt}V(p(t), x(t)) = \frac{-x(t)^2 + x(t)u(t) +x(t)\frac{(-x(\hat t-\theta) + u(\hat t))\theta}{1+p_0^2}}{1+p_0^2}\\[3mm]
& &\qquad \le
\frac{-x(t)^2 + \frac14 x(t)^2 + u(t)^2 + \abs{x(t)}\,(\max_{\tau\in [0, t]}\normo{x_\tau} + \norm{u})\theta}{1+p^2_0}\\[3mm]
& &\qquad\le
\frac{-\frac34x(t)^2 + u(t)^2 + \theta\cdot(\max_{\tau\in [0, t]}\normo{x_\tau}^2 + \frac14x(t)^2 + \norm{u}^2)}{1+p_0^2}\\
& &\qquad\le 
\frac{-\frac{1}2x(t)^2 + (1+\theta)\norm{u}^2 + \theta\cdot
\max_{\tau\in [0, t]}\normo{x_\tau}^2}{1+p_0^2}.
\eeqn
Hence 
\[
\frac{d}{dt}V(p(t), x(t))\le\frac{-V(p(t), x(t)) + (1+\theta)\norm{u}^2 + \theta\cdot\max_{\tau\in[0, t]}\normo{x_\tau}^2}{1+p_0^2}.
\]
By a standard comparison method with solutions of the linear system $\dot z(t) = -\frac{z(t)}{1+p_0^2} + v(t)$ with $v(t)\ge 0$,
one can conclude that
\beqn
V(p(t), x(t)) &\hspace{-4pt}\le&\hspace{-7pt} V(p_0, x(0))e^{-\frac{t}{1+p_0^2}}
+\theta
\max_{\tau\in[-\theta, t]}\normo{x_\tau}^2
+ (1+\theta)\norm{u}^2\\[3mm]
{}\hspace{-6pt}&\le&\hspace{-7pt}
V(p_0, x(0))e^{-\frac{t}{1+p_0^2}}
+
2\theta\!\!\max_{\tau\in[-\theta, t]}V(p(\tau), x(\tau))
+ (1+\theta)\norm{u}^2.
\eeqn
Applying Statement 2 of Proposition \ref{p-raz-ios}, one concludes that if $\theta<\frac12$, then the system with the new (delay-free) output map $V(p, x)$ is \olios, that is, for some $\beta\in\kl$ and some $\rho, \gamma\in\Nn$, it holds that
\[
V(p(t), x(t))\le \beta\left(\normo{x_0}, \ \frac{t}{1+\rho(\normo{(x_0, p_0)})}\right)
+ \gamma(\norm{u}), \quad\forall\,t\ge 0.
\]
Consequently, 
\[
\abs{x(t)}
\le \left(2\beta\left(\normo{x_0}, \ \frac{t}{1+\rho(\normo{(x_0, p_0)})}\right)
+ 2\gamma(\norm{u})\right)^{1/2}, \quad\forall\,t\ge 0.
\]
In terms of the original system \rref{e-ex2}, this means that the system is \iss\ but the decay rate of $x(t)$ is affected by the value of the parameter $p$. 
\eex

\section{IOS and OL-IOS}\label{s-4}

As remarked earlier, the major difference between \ios\ and \olios\ lies on whether
the overshoots of output signals are dominated by $\abs{y(0)}$ or by
$\normo{x_0}$. The stricter \olios\ condition leads to a significant step in the
converse Lyapunov theorems for a variety of output stability properties for systems
without delays.  As a preparation to develop the converse Lyapunov theorems for delay systems,
we present the next result on how \ios\ and \olios\ are related to systems with delays.

For continuous function $q:\xc\rightarrow\R$, define the set $\Zz_q\subseteq\xc$ by
\[
\Zz_q = \{\xi\in\xc: \ q(\xi)=0\}.
\]
\bd{d-redefinition}
Let $h$ be the output map of \rref{e-syso}.
A continuous map $\bar h: \xc \rightarrow \R_{\geq 0}$ is said to be an
output redefinition of $h$ if there exists some $\ab, \at\in\ki$ such that
\[
\underline{\alpha} \left(\abs{h(\xi)} \right) \leq
  \enorm{\bar h(\xi)}\leq  \overline{\alpha} \left( \normx{\xi}\right), \quad
  \forall\xi\in\xc.
\]
\eds


\comment{
\bd{olios-redefinition}
A system as in \rref{e-syso} is \olios\ {\it under output redefinition}
if for some output redefinition $h_0$ of $h$, the system
\be{e-syso-redef}
\dot x(t)= f(x_t, u(t)), \quad y(t) = h_0(x_t)
\ee
is \olios.
\eds
}

A redefinition $\bar h$ of $h$ is roughly a map that dominates $h$. The result below states that the \ios\ property implies the \olios\ with an
output redefinition.

\comment{
\rcomment{By combining Theorem \ref{olios-under-or} with \rref{e-implication2}, the interconnections of the output stability notions can be summarized as below: }
\begin{center}
\begin{tikzpicture}[
     implies/.style={double,double equal sign distance,-implies},
     dot/.style={shape=circle,fill=black,minimum size=2pt,
                 inner sep=0pt,outer sep=2pt},
]

    \node (siios) [box] {\siios};
    \node (olios) [box, right of=siios, xshift=3.5 cm] {\olios};
    \node (ios) [box, right of=olios, xshift=3.5 cm] {\ios};
    \node (ros) [box, right of=ios, xshift=3.5 cm] {\ros};
    \draw (ios) edge[implies] node[above] {Lemma \ref{l-ios-ros}} (ros);
    \draw (siios) edge[implies] node[above] {Definition \ref{d-ios}} (olios);
    \draw (olios) edge[implies] node[above] {Definition \ref{d-ios}} (ios);
    \draw (ios) edge[implies, bend left] node[below] {Theorem \ref{olios-under-or}} (olios);
\end{tikzpicture}
\end{center}
}

\bt{olios-under-or}
The following are equivalent for a \rfc\ system \rref{e-syso}:
\begin{enumerate}
\item[(a)] The system  is \ios.
\item[(b)] The system  is \olios\ under a output redefinition.
\end{enumerate}
Furthermore, if $h$ in \rref{e-sysoo} is locally Lipschitz, then the output redefinition $\bar h$ can be chosen to be locally Lipschitz on $\xc\setminus\Zz_{\bar h}$.
\ets

The proof of the theorem will be given in Section \ref{s-pf1}.

\bex{ex-red-olios} To illustrate the meaning of Theorem \ref{olios-under-or}, consider the system:
\be{e-ex1}
\begin{array}{ccl}
\dot x_1(t) &=& -x_1(t) + x_1(t-\theta),\\[2mm]
\dot x_2(t) &=& -\frac{x_2(t)}{1+\abs{x_1(t)}^2} - x_3(t)\normo{x_t}^2,\\[2mm]
\dot x_3(t) &=& - \frac{x_3(t)}{1+\abs{x_1(t)}^2}\, + \displaystyle{x_2(t)\normo{x_t}^2} + \frac{u(t)}{1+\abs{x_1(t)}^2},
\end{array}
\ee
with $y(t) = x_3(t)$ as output, i.e., $h(\xi) = \xi_3(0)$. Note that 
\[
\normo{x_t} = \max_{s\in [-\theta, 0]}\abs{(x_1(s), x_2(s), x_3(s))}.
\]
Let
\[
W_1(\xi_1) = \frac12 \xi_1(0)^2 + \frac12\int_{-\theta}^0\xi_1^2(s)\,ds.
\]
Note that
\[
\frac{1}{2}|\xi_1(0)|^2\le W_1(\xi_1)\le\frac{1+\theta}{2}\normo{\xi_1}^2.
\]
Along the trajectories of the $x_1$-subsystem,
\[
W_1((x_1)_t) = \frac12x_1(t)^2 + \frac12\int_{t-\theta}^tx_1(s)^2\,ds,
\]
and
\beqn
& &\frac{d}{dt}W_1((x_1)_t) = -x_1(t)^2 + x_1(t)x_1(t-\theta) + \frac12x_1(t)^2 - \frac12x_1(t-\theta)^2\\
& &\le {} -x_1(t)^2 + \frac12x_1(t)^2 + \frac12x_1(t-\theta)^2 + \frac12x_1(t)^2 - \frac12x_1(t-\theta)^2\le 0.
\eeqn
It follows that $W_1((x_1)_t)\le W_1(\xi_1)$ for all $t\ge 0$, and hence,
\[
\frac{\abs{x_1(t)}^2}{2} \le W_1((x_1)_t)
\le \frac{1+\theta}{2}\,\normo{\xi_1}^2\,, \quad\forall\, t\ge 0.
\]
Thus, $\abs{x_1(t)}\le \sqrt{1+\theta}\,\,\normo{\xi_1}$ for all $t\ge 0$ along every trajectory of the $x_1$-subsystem with the initial state $\xi_1$. Let
$W_2(\xi) = \frac{1}{2}(\xi_2(0)^2 + \xi_3(0)^2)$.   Along the trajectory of the system, it holds that
\[
W_2(x_t) =\frac12x_2(t)^2 + \frac12x_3(t)^2,
\]
and
\beqn
\frac{d}{dt}W_2(x_t)) &=&-\frac{x_2(t)^2}{1+\abs{x_1(t)}^2} - \frac{x_3(t)^2}{1+\abs{x_1(t)}^2} +  \frac{x_3(t)u(t)}{1+\abs{x_1(t)}^2} \\
&\le& -\frac{x_2(t)^2}{1+\abs{x_1(t)}^2} - \frac{x_3(t)^2}{1+\abs{x_1(t)}^2} + 
\frac{\frac12(x_3(t)^2 + u(t)^2)}{1+\abs{x_1(t)}^2}\\
&\le&\frac{1}{1+\abs{x_1(t)}^2}\left(-\frac12x_2(t)^2 - \frac12x_3(t)^2 + \frac12u(t)^2\right)\\
&\le&\frac{1}{1+\abs{x_1(t)}^2}\left(-W_2(x_t) + \frac12u(t)^2\right).
\eeqn
It then follows that
\beqn
W_2(x_t)\ge \abs{u}^2\quad \Rightarrow\
\frac{d}{dt}W_2(x_t)&\le& -\frac{1/2}{1+\abs{(x_1)_t}^2}W_2(x_t)\\
&\le& -\frac{1/2}{1+c\normx{\xi_1}^2}W_2(x_t),
\eeqn
where $c=\sqrt{1+\theta}$.
The following estimate can then be derived:
\[
W_2(x_t)\le W_2(\xi)e^{-\frac{t}{2(1+c\normo{\xi_1}^2)}} + \norm{u}^2.
\]
With the inequality that $\abs{x_3(t)}^2\le2W_2(x_t)$, one concludes that the system is \ios\ with $x_3$ as output.  However, notice that the system is not \olios\ as $y(t)\not\equiv0$ even if $y_0 =\xi_2=0$ and $u\equiv0$.  On the hand, it should not be hard to see that the system is \olios\ under the output redefinition with $W_2(\xi)$ as output.

This example can also be used to illustrate the meaning of \olios: with $W_2(\xi)$ as the redefined output, the overshoot of the output function is dominated by $W_2(\xi)$ and $\norm{u}$; however, large initial values of  $x_1(t)$ can slow down the decay rate of the output map.
\eex

\section{Robust Output Stability Under Output Feedback}\label{s-5}

If a system as in \rref{e-syso} is \ios, it is zero-input output stable.
In this section, we provide an extension of a result given
in \cite{sw-ios-scl} to delay systems to show that the \ios\ property allows a $\ki$
stability margin under output feedback.

\subsection{Robust Global Asymptotic Output Stability}

Let $\Uu$ be a subset of $\R^m$.  Consider a system as follows:
 \begin{subequations}\label{e-sysro}
\begin{align}
\dot x(t) &= f_0(x_t, d(t)), \label{e-sysrox}\\
y(t) &= h(x_t), \label{e-sysroo}
\end{align}
\end{subequations}
where for each $t$, $x(t)\in\R^n, \ d(t)\in \Uu $, and $y(t)\in\R^p$.
The map $f_0: \xo^n\times \Uu \rightarrow \R^n$ is locally Lipschitz, and $h:\xc\rightarrow\R^p$ is continuous, satisfying the \olgs\ property \rref{bound-h}. We will use $x(t, \xi, d)$ to denote the trajectory of \rref{e-sysro} corresponding to the initial state $\xi\in\xc$ and $d\in\MmU$; and correspondingly, $y(t,\xi, d):=h(x(t, \xi, d))$.  The system is forward complete if all trajectories of the system are defined on $[0, \infty)$.

\bd{d-rgaos}
Consider a forward complete system \rref{e-sysro}. The system is said to be
\begin{enumerate}
\item
{\it robustly globally asymptotically output stable} (\rgaos) if there exists $\beta\in\kl$ such that
\be{e-rgaos}
\abs{y(t, \xi, d)}\le \beta(\normo{\xi}, t)\qquad\forall\, t\ge 0, \ \forall\,\xi\in\xc, \ \forall\,d\in\MmU;
\ee
\item
{\it output-Lagrange \rgaos} (\olrgaos) if the system is \rgaos, and additionally, there exists some $\sigma\in\kk$ such that
\be{e-olrgaos}
\abs{y(t, \xi, d)}\le\sigma(\abs{h(\xi)}); \qquad \forall\, t\ge 0, \ \forall\,\xi\in\xc, \ \forall\,d\in\MmU;
\ee
\item
{\it state-independent \rgaos} (\sirgaos) if there exists $\beta\in\kl$ such that
\be{e-sirgaos}
\abs{y(t, \xi, d)}\le \beta(\abs{h(\xi)}, t)\qquad\forall\, t\ge 0, \ \forall\,\xi\in\xc, \ \forall\,d\in\MmU;
\ee
\end{enumerate}
\eds
As in the case of \olios, system \rref{e-sysro} is \ol-\rgaos\ if and only if there exist $\beta\in\kl$ and $\rho\in\kk$ such that
\[
\abs{y(t, \xi, d)}\le \beta\left(\abs{h(\xi)}, \ \frac{t}{1+\rho(\normo{\xi})}\right),\quad\forall\,t\ge 0, \ \forall\,\xi\in\xc, \
\forall\,d\in\MmU.
\]

As in the case of \olios\ and \siios, when the output map $h$ is delay-free, i.e., when $h(\xi)=h_0(\xi(0))$, then we say that system \rref{e-sysro} is \olrgaos\ in the history norm if
the system is \rgaos\ and satisfies the following \ol-property with $H$ as the output, where $\sigma\in\kk$:
\[
\abs{y(t, \xi, d)}\le \sigma\left(H(\xi)\right),\qquad\forall\,t\ge 0, \ \forall\,\xi\in\xc, \ \forall\,d\in\MmU,
\]
(recall that $H(\xi) = \max_{s\in [-\theta, 0]}\abs{h_0(\xi(s))}$), and system \rref{e-sysro} is \sirgaos\ if there exists some $\beta\in\kl$ such that
\[
\abs{y(t, \xi, d)}\le \beta(H(\xi),\ t), \qquad\forall\,t\ge 0, \ \forall\,\xi\in\xc, \ \forall\,d\in\MmU.
\]

When the set $\Uu$ is bounded, the \rfc\ property \rref{e-rfc} for \rref{e-sysro} means that there exist $\chi_1, \chi_2\in\kk$ and $c>0$ such that
\[
\abs{x(t, \xi, d)}\le \chi_1(t) + \chi_2(\normo{\xi}) + c, \quad\forall\,t\ge 0, \ \forall\,\xi\in\xc, \ \forall\,d\in\MmU.
\]

\subsection{Stability Margin Under Output Feedback}
To a given system as in \rref{e-syso} and a locally Lipschitz function
$\lambda:[0, \infty)\rightarrow [0, \infty)$, consider
the system with disturbances $d(\cdot)\in\MO$:
\beq\label{alt-sys}
 \dot{x} (t) &=& g\left(x_t, d(t) \right)
:= f\left(x_t, d(t) \lambda (\abs{y(t)}) \right),\\
y &=&h(x_t), \nonumber
\eeq
where $\Omega$ denotes the closed unit ball of $\R^m$.
We use $x_\lambda(\cdot, \xi, d)$ and $y_\lambda(\cdot,
 \xi, d)$ to denote correspondingly the trajectory and the output of \rref{alt-sys}.
Note that \rref{alt-sys} needs not be forward complete
even if the original system \rref{e-syso} is. On the other hand,
on the maximum interval $[0, T_{\xi,d})$ of a trajectory of
\rref{alt-sys}, it holds that $x_\lambda(t,\xi, d) =
x(t, \xi, u)$, the trajectory of \rref{e-syso} with the initial state
 $\xi$ and the input $u$ given by $u(t) =
 d(t)\lambda(\abs{y_\lambda(t, \xi,  d)})$.  As discussed earlier, if a system is \ios, then the corresponding zero-input system is \gaos.  The next result shows that the global asymptotic output stability is robust under a $\ki$ margin of output feedback.

\bp{p-robust}
Consider a system as in \rref{e-syso} and assume that it is \rfc.  Further assume that the output map $h$ is locally Lipschitz. 
\begin{enumerate}
\item If the system is \ios, then there exists $\lambda\in\ki$ which is smooth such that the corresponding system \rref{alt-sys} is \rfc\ and \rgaos.
\item If the system is \olios, then there exists $\lambda\in\ki$ which is smooth such that the corresponding system \rref{alt-sys} is \rfc\ and \olrgaos.
\item If the system is \siios, then there exists $\lambda\in\ki$ which is smooth such that the corresponding system \rref{alt-sys} is \rfc\ and \sirgaos.\mybox
\end{enumerate}
\eps

\comment{
\bp{l-ios-ros}
Consider a system as in \rref{e-syso} and assume that it is \rfc.  Assume that the output map $h$ is locally Lipschitz. If it is
\ios, then there exists $\lambda\in\ki$ which is locally Lipschitz such that the corresponding system \rref{alt-sys} is \rfc\ and \rgaos.
\ep

As remarked earlier, the \olios\ property turns out to be the most critical one
for establishing the Lyapunov descriptions for both \ios\ and
\olios\ (see \cite{klw-ifac17}). The proof of Proposition \ref{l-ios-ros} also
requires the next result which is by itself interesting.

\bp{l-olios-olros}
Consider a \rfc\ system as in \rref{e-syso}. Assume that the output map $h$ is locally Lipschitz.
If it is \olios\ (\siios\ respectively), then there exist a smooth
$\lambda\in\ki$ such that the corresponding system \rref{alt-sys} is \rfc\ and \olrgaos\ (\sirgaos\ respectively).
\comment{
and additionally, satisfies the following robust \ol\ property:
\[
\abs{y_{\lambda}(t, \xi, d)}
\leq \sigma\left(\abs{h(\xi)}\right)
\]
for all $t\ge 0$, $\xi\in\xc$, and $d\in\MO$.
}
\ep
}

The proof of Proposition \ref{p-robust} will be given in Section \ref{s-6}.

\section{Proofs}\label{s-6}

In this section we provide proofs of the results given in Sections \ref{s-raz}, \ref{s-4}, and \ref{s-5}.

\subsection{Proof of \protect{Proposition \ref{p-raz-ios}}}

To prove the proposition, we first prove the following lemma 
for a forward complete system \rref{e-syso} where $h(\xi) = h_0(\xi(0)$ with $h_0(0)=0$. Recall that the associated map of output history $H:\xc\rightarrow\R_{\ge 0}$ is defined by $H(\xi) = \max_{\tau\in [-\theta, 0]}\abs{h_0(\xi(s))}$, that is, $H(\xi) = \normo{y_0}$.
\bl{l-raz-new1}
Let $\alpha:\xc\rightarrow\R_{\ge 0}$ be a continuous function with $\alpha(0)=0$.
Assume that there exist $\beta\in\kl$, $\rho, \kappa, \gamma\in\Nn$ such that
\beq
\abs{y(t)} &\le& \max\left\{\beta\left(\alpha(\xi), \ \frac{t}{1+\rho(\normo{\xi})}\right),\right.\nonumber\\ 
& &\qquad\qquad\left. \kappa\left(\max_{\tau\in [-\theta, t]}\abs{y(\tau)}\right),\ \gamma(\norm{u})\right\},\quad\forall\,t\ge 0.\label{e-raz5}
\eeq
Then there exists some $\beta_1\in\kl$ such that
the following holds along all trajectories:
\beq
\normo{y_t^{\xi, u}} &\le& 
\max\left\{\beta_1\left(\tilde\alpha(\xi), \ \frac{t}{1+\rho(\normo{\xi})}\right), \right. \nonumber\\ 
& &\qquad\qquad\left.\kappa\left(\max_{\tau\in [0, t]}\normo{y_\tau} \right),\ \gamma(\norm{u})\right\}, \quad\forall\,t\ge 0,
\label{e-ios-norm}
\eeq 
where $\tilde\alpha(\xi) = \max\{\alpha(\xi), \ H(\xi)\}$.
\el

\bpr
Assume that \rref{e-raz5} holds along all trajectories. Without loss of generality, we may assume that $\beta(r, 0)\ge r$ (otherwise replace $\beta(r, t)$ by $\beta(r, t) + \frac{r}{1+t}$.)

Let $\xi\in\xc$ and $u\in\Mm$ be given, and we use $y(t)$ to denote the corresponding output function.  It can be seen that if $t\ge\theta$, then 
\[
\normo{y_t^{\xi, u}}\le
\max\left\{\beta\left(\alpha(\xi), \ \frac{t-\theta}{1+\rho(\normo{\xi})}\right),\ \kappa\left(\max_{\tau\in [-\theta, t]}\abs{y(\tau)} \right),\ \gamma(\norm{u})\right\}.
\]
For $t\in [0, \theta)$,
\beqn
\normo{y_t^{\xi, u}} &\le&
\max\left\{H(\xi), \ 
\max\left\{\beta(\alpha(\xi), 0), \kappa\left(\max_{\tau\in[-\theta, t)}\abs{y(\tau)}\right),
\gamma(\norm{u})\right\}\right\}\\
&\le&
\max\left\{\beta(\tilde\alpha(\xi), 0), \ \kappa\left(\max_{\tau\in [-\theta, t]}\abs{y(\tau)} \right),\ \gamma(\norm{u})\right\}.
\eeqn
\comment{
Assume $t\in[0, \theta)$ and $\theta\in[-\theta, 0]$.  If $t-s\ge 0$, then
\beqn
\abs{y(t-s)}&\le& \max\left\{\beta(\normo{\xi}, \ t-s), \ \kappa\left(\max_{\tau\in [-\theta, t]}\abs{y(\tau)}) \right),\ \gamma(\norm{u})\right\}\\
&\le&
\max\left\{\beta(\normo{\xi}, \ 0), \ \kappa\left(\max_{\tau\in [-\theta, t]}\abs{y(\tau)}) \right),\ \gamma(\norm{u})\right\};
\eeqn
and for $t-s<0$, it holds that
\[
\abs{y(t-s)} = \abs{h_0(x(t-s)))}\le \pi(x(t-s)) \le\pi(\normo{\xi}),
\]
where $\pi$ is as in \rref{bound-h}. Hence, for all $0\le t < \theta$,
\[
\normo{y_t^{\xi, u}}\le \max\left\{\pi_1(\normo{\xi}), \ \kappa\left(\max_{\tau\in [-\theta, t]}\abs{y(\tau)}) \right),\ \gamma(\norm{u})\right\}, 
\]
where $\pi_1(r) = \max\{\beta(r, 0), \pi(r)\}$.
}
Applying Proposition 7 in \cite{eds-iISS-scl98} to the $\kl$-function $\beta$, one sees that there exist some $p\in\ki$ and $q\in\ki$ such that $p(\beta(r, t))\le q(r)e^{-t}$, and thus, for any $r\ge 0$ and $s\ge 0$,
\[
p\left(\beta\left(r, \ \frac{t-\theta}{1+\rho(s)}\right)\right)\le e^{\frac{\theta}{1+\rho(s)}}q(r)e^{-\frac{t}{1+\rho(s)}}
\le e^{\theta}q(r)e^{-\frac{t}{1+\rho(s)}},\quad\forall\,t\ge \theta.
\]
Let
\[
\beta_1(r, t) = p^{-1}\left(e^{\theta}q(r)e^{-t}\right).
\]
Then the following holds along all trajectories:
\[
\normo{y_t}\le
\max\left\{\beta_1\left(\tilde\alpha(\xi), \ \frac{t}{1+\rho(\normo{\xi})}\right), \ \kappa\left(\max_{\tau\in [-\theta, t]}\abs{y(\tau)} \right),\ \gamma(\norm{u})\right\}, 
\]
for all $t\ge 0$.
Note that this is equivalent to
\[
\normo{y_t}\le
\max\left\{\beta_1\left(\tilde\alpha(\xi), \ \frac{t}{1+\rho(\normo{\xi})}\right), \ \kappa\left(\max_{\tau\in [0, t]}\normo{y_\tau} \right),\ \gamma(\norm{u})\right\},
\]
for all $t\ge 0$, all $\xi\in\xo$, and all $u\in\Mm$.
\epr

\noindent{\it Proof of Proposition \ref{p-raz-ios}, Statement 1.}
Suppose \rref{e-raz-ios} holds for some $\beta\in\kl$, $\kappa, \gamma\in\Nn$ with $\kappa(s) < s$ for all $s>0$. 
By Lemma \ref{l-raz-new1}, there is some $\beta_1\in\kl$ such that
\beqn
\normo{y_t} &\le&
\max\left\{\beta_1\left(\tilde\alpha(\normo{\xi}), \ t\right), \ \kappa\left(\max_{\tau\in [0, t]}\normo{y_\tau} \right),\ \gamma(\norm{u})\right\}\\
&\le&
\max\left\{\beta_2\left(\normo{\xi}, \ t\right), \ \kappa\left(\max_{\tau\in [0, t]}\normo{y_\tau} \right),\ \gamma(\norm{u})\right\},
\eeqn
for all $t\ge 0$, all $\xi\in\xo$, and all $u\in\Mm$, where $\tilde\alpha(r) = \max\{r, \ \pi(r)\}$, $\pi$ is as in \rref{bound-h},  and $\beta_2(r, t) = \beta_1(\tilde\alpha(r), t)$. It then follows that
\[
\normo{y_t}\le\max\left\{
\beta_2\left(\norm{x_{t/2}}, \ \frac{t}2\right),\
\kappa\left(\max_{\tau\in[t/2, t]}\normo{y_\tau}\right),\ \gamma(\norm{u})\right\},
\]
for all $t\ge 0$.  Since system \rref{e-sysox} is assumed to satisfy the \gs\ condition, one sees that for some $\sigma\in\ki$ and $\mu\in\Nn$,
\[
\normo{x_\tau}\le\max\left\{\sigma(\normo{\xi}), \ \mu(\norm{u})\right\}, \quad\forall\,\tau\ge 0.
\]
It follows that
\beqn
\normo{y_t} &\le& \max\left\{
\beta_2\left(\sigma(\normo{\xi}), \frac{t}2\right),\
\beta_2(\mu(\norm{u}), \ 0),\right. \\
& &\qquad\qquad \qquad\left. \kappa\left(\max_{\tau\in[t/2, t]}\normo{y_\tau}\right),\ \gamma(\norm{u})\right\}\\
&\le& \max\left\{
\beta_3(\normo{\xi}, t),\  \kappa\left(\max_{\tau\in[t/2, t]}\normo{y_\tau}\right),\ \gamma_1(\norm{u})\right\}\\
\eeqn
for all $t\ge 0$, where $\beta_3(r, t) = \beta_2(\sigma(r), t/2)$ and $\gamma_1(r) = \max\{\gamma(r), \mu(r)\}$.
Applying Lemma 5.4 in \cite{twj-dcdis} (with $a=\normo{\xi}$, $\mu=\frac12$, and $b=\gamma_1(\norm{u})$), one sees that there exists some $\beta_4\in\kl$ which depends only on $\beta_3$ such that
\[
\normo{y}\le \max\left\{\beta_4(\normo{\xi}, \ t), \ \gamma_1(\norm{u})\right\},\quad\forall\,t\ge 0.
\]
This proves Statement 1 of Proposition \ref{p-raz-ios}.

\bigskip

\noindent{\it Proof of Proposition \ref{p-raz-ios}, Statement 2.}
Assume \rref{e-raz-ol} holds for some $\beta\in\kl$, $\rho, \kappa, \gamma\in\Nn$.  By Lemma \ref{l-raz-new1},
there exists some $\beta_5\in\kl$ such that the following holds along all trajectories:
\be{e-yt}
\normo{y_t}\le\max\left\{
\beta_5\left(H(\xi), \ \frac{t}{1+\rho(\normo{\xi})}\right),\ 
\kappa\left(\max_{\tau\in [0, t]}\normo{y_\tau} \right),\ \gamma(\norm{u})\right\},
\ee
for all $t\ge 0$.

Let $T>0$ be given.  Then
\[
\normo{y_t}\le\max\left\{
\beta_5\left(H(\xi), \ 0\right),\ 
\kappa\left(\max_{\tau\in [0, T]}\normo{y_\tau} \right),\ \gamma(\norm{u})\right\},
\]
for all $t\in[0, T]$. Hence,
\[
\max_{\tau\in[0, T]}\normo{y_\tau}\le\max\left\{
\beta_5\left(H(\xi), \ 0\right),\ 
\kappa\left(\max_{\tau\in [0, T]}\normo{y_\tau} \right),\ \gamma(\norm{u})\right\}.
\]
Since $\kappa(s)<s$ for all $s>0$, it follows that
\[
\max_{\tau\in[0, T]}\normo{y_\tau}\le \max\left\{
\beta_5\left(H(\xi), \ 0\right),\ \gamma(\norm{u})\right\}, 
\]
and in particular,
\be{e-H3}
\normo{y_T}\le \max\left\{
\beta_5\left(H(\xi), \ 0\right),\ \gamma(\norm{u})\right\}. 
\ee
Observing that \rref{e-H3} for all $T\ge 0$, one gets the following \olgs\ property in the history norm:
\be{e-Hol}
\normo{H(x_t)}\le \max\left\{\sigma_1(H(\xi)), \ \gamma(\norm{u})\right\}, \quad\forall\,t\ge 0,
\ee
where $\sigma_1(r) = \beta_5(r, 0)$.

Let $\xi\in\xc$ and $u\in\Mm$.  Applying again the semi-group property to \rref{e-yt}, one sees that
\beqn
\normo{y_t} &\le& \max\left\{
\beta_5\left(\normo{y_{t/2}}, \ \frac{\frac{t}2}{1+\rho(\normo{x_{t/2}})}\right),\right.\\
& &\hspace{1.2in}\left.
\kappa\left(\max_{\tau\in[t/2, t]}\normo{y_\tau}\right),\ \gamma(\norm{u})\right\},
\eeqn
for all $t\ge 0$.  By \rref{e-Hol}, it follows that
\beq
\normo{y_t} &\le& \max\left\{
\beta_5\left(\sigma_1(H(\xi)), \ \frac{\frac{t}2}{1+\rho(\normo{x_{t/2}})}\right),\right.\nonumber\\
& &\hspace{1.2in}\left.
\kappa\left(\max_{\tau\in[t/2, t]}\normo{y_\tau}\right),\ \gamma_2(\norm{u})\right\},\label{e-yt2}
\eeq
for all $t \ge 0$, where $\gamma_2(r) = \max\{\gamma(r), \beta_5(\gamma(r), 0)\}$.

By the \ubibs\ condition imposed on system \rref{e-syso}, one sees that for some $\sigma_2\in\ki$, some
$\mu_2\in\Nn$ and some $c\ge 0$, it holds that
\[
\normo{x_\tau}\le\max\{\sigma_2(\normo{\xi}) + c, \ \mu_2(\norm{u})\}, \quad\forall\tau\ge 0.
\]
It follows from \rref{e-yt2} that
\beq
& &\normo{y_t} \le \max\left\{
\beta_5\left(\sigma_1(H(\xi)), \ \frac{\frac{t}2}{1+\max\{\rho_3(\normo{\xi}), \mu_3(\norm{u})\}} \right),\right.\nonumber\\[3mm]
& &\hspace{1.7in}\left.
\kappa\left(\max_{\tau\in[t/2, t]}\normo{y_\tau}\right),\ \gamma_2(\norm{u})\right\},\nonumber\\[3pt]
& &\le
\max\left\{\beta_5\left(\sigma_1(H(\xi)), \ \frac{\frac{t}2}{1+\rho_3(\normo{\xi})}\right),
\beta_5\left(\sigma_1(H(\xi)), \ \frac{\frac{t}{2}}{1+\mu_3(\norm{u})}\right),\right. \nonumber\\[3pt]
& &\hspace{1.7in}\left.
\kappa\left(\max_{\tau\in[t/2, t]}\normo{y_\tau}\right),\ \gamma_2(\norm{u})\right\},\label{e-raz7}
\eeq
where $\rho_3(r) = \rho(\sigma_2(r)+c)$, and $\mu_3(r) = \rho(\mu_2(r))$.
For the first term in the above expression, it holds that
\beqn
\beta_5\left(\sigma_1(H(\xi)), \ \frac{\frac{t}2}{1+\rho_3(\normo{\xi})}\right)
&\le&
\beta_5\left(\sigma_1(\pi(\normo{\xi})), \ \frac{\frac{t}2}{1+\rho_3(\normo{\xi})}\right)\\
&=&
\hat\beta_5(\normo{\xi}, t),
\eeqn
where $\hat \beta_5(r, t) =\beta_5(\sigma_1(\pi(r)), \frac{t/2}{1+\rho_3(r)})$.  Note that $\hat \beta_5\in\kl$.
For the second term in the above expression, by considering the two cases when $H(\xi)\ge \norm{u}$ and when $H(\xi)\le \norm{u}$, one gets the following:
\beqn
& &\beta_5\left(\sigma_1(H(\xi)), \ \frac{\frac{t}2}{1+\mu_3(\norm{u})}\right)\\
& &{}\le
\max\left\{
\beta_5\left(\sigma_1(H(\xi)), \ \frac{\frac{t}2}{1+\mu_3(H(\xi))}\right), \ 
\beta_5\left(\sigma_1(\norm{u}), \ 0\right),
\right\}\\
& &{}\le\max\left\{
\beta_5\left(\sigma_1(\pi(\normo{\xi})), \ \frac{\frac{t}2}{1+\mu_3(\pi(\normo{\xi}))}\right),
\ \beta_5(\sigma_1(\norm{u}), 0)\right\}\\
& & {}\le\max\left\{
\tilde\beta_5(\normo{\xi}, \ t),\ \tilde\mu_3(\norm{u})\right\},
\eeqn
where $\tilde\beta_5(r, t) = \beta_5(\sigma_1(\pi(r)), \ \frac{t/2}{1+\mu_3(\pi(r))})$, and
$\tilde\mu_3(r)=\beta_5(\sigma_1(r), 0)$.
Thus, \rref{e-raz7} leads to the following:
\be{e-raz8}
\normo{y_t} \le \max\left\{\beta_6(\normo{\xi}, t), \ 
\kappa\left(\max_{\tau\in[t/2, t]}\normo{y_\tau}\right),\ \gamma_3(\norm{u})\right\},
\ee
for all $t\ge 0$, where
\[
\beta_6(r, t)=\max\{\hat\beta_5(r, t), \tilde\beta_5(r, t)\},\quad 
\gamma_3(r) = \max\{\gamma_2(r), \tilde\mu_3(r)\}. 
\]
Again, by applying Lemma 5.4 in \cite{twj-dcdis}, one concludes that there exists some $\beta_7\in\kl$ such that 
\be{e-raz9}
\normo{y_t} \le \max\left\{\beta_7(\normo{\xi}, t), \ 
 \gamma_3(\norm{u})\right\},
\ee
for all $t\ge 0$.  Since $\xi\in\xc$ and $u\in\Mm$ were chosen arbitrarily, one concludes that the \ios\ property holds for the system.  Combining this with the \olgs\  property \rref{e-Hol}, one sees that the system is \olios\ in the history norm.

\medskip 

\noindent{\it Proof of Proposition \ref{p-raz-ios}, Statement 3.}
The proof of Statement 3 follows the same steps as in the proof of Statement 1. Assume 
\rref{e-raz-siios} holds for some $\beta\in\kl$, $\kappa, \gamma\in\Nn$ with $\kappa(s)<s$ for all $s>0$. Then there exists some $\beta_1\in\kl$ such that the following holds along all trajectories:
\[
\normo{y_t}\le \max\left\{\beta_1(\normo{y_0}, \ t), \kappa\left(\max_{\tau\in [0, t]}\normo{y_\tau} \right),\ \gamma(\norm{u})\right\},
\] 
for all $t\ge 0$.  Let $\xi\in\xo$ and $u\in\Mm$ be given. Then
\be{e-raz10}
\normo{y_t} \le \max\left\{\beta_1(\normo{y_{t/2}}, \ t/2), \kappa\left(\max_{\tau\in [t/2, t]}\normo{y_\tau} \right),\ \gamma(\norm{u})\right\},
\ee
for all $t\ge 0$.
With $\rho(r)=0$ for all $r\ge0$, one still gets the \olgs\ property \rref{e-Hol} (whose proof does not use the \ubibs\ condition), and thus,
\beqn
\normo{y_t} &\le& \max\left\{\beta_1(\sigma_1(\normo{y_0}), \ t/2), \kappa\left(\max_{\tau\in [t/2, t]}\normo{y_\tau} \right),\ \gamma_4(\norm{u})\right\}\\
&\le&  \max\left\{\hat\beta_1(H(\xi), \ t), \ \kappa\left(\max_{\tau\in [t/2, t]}\normo{y_\tau} \right),\ \gamma_4(\norm{u})\right\}, \quad\forall\,t\ge 0,
\eeqn 
where $\gamma_4(r) = \max\{\gamma(r), \beta_1(\gamma(r), 0)\}$, and $\hat\beta_1(r, t)=\beta_1(\sigma_1(r), t/2)$.
From this one concludes that for some $\beta_8\in\kl$, the following \siios\ estimate holds along all trajectories:
\[
\normo{y_t} \le \max\left\{\beta_8(H(\xi), \ t), \ \gamma_4(\norm{u})\right\}, \quad\forall\,t\ge 0.
\]
The proof is completed by noticing that $\xi\in\xo, u\in\Mm$ can be chosen arbitrarily.
\qed

\subsection{Proof of \protect{Theorem \ref{olios-under-or}}}\label{s-pf1}

The implication (b) $\Rightarrow$ (a) follows immediately from the condition
$\abs{h(\xi)}\le \ab^{-1}(\abs{h_0(\xi)})$.
This leaves the implication  (a) $\Rightarrow$ (b) as the main
component of the proof, which we present below.

Suppose that system \rref{e-syso} is \ios, and thus, there
exist $\beta
\in \kl$ and $\gamma \in \ki$ such that
\be{e-pf1}
\abs{y(t, \xi, u)} \leq \max \{ \beta
\left( \normx{\xi}, t \right), \gamma (\norm{u}) \}\quad
\forall\,t\ge 0
\ee
along all trajectories. Define $\bar{h} : \xc \rightarrow \R $ by
\[
\bar{h} (\xi ) := \sup_{t \geq 0, u \in \Mm} \left\{
\max \{ \abs{y(t, \xi, u)} -\gamma (\norm{u}), 0 \} \right\}.
\]
Note that with $u\equiv 0$,
\be{e-pf10}
\abs{ y(0, \xi, u) } - \gamma (\norm{u}) = \abs{h(\xi)}\ge 0\quad
\forall\,\xi.
\ee
Hence, the map $\bar{h}$ is equivalently defined by
\be{1st-or}
\bar{h} (\xi ) := \sup_{t \geq 0, u \in \Mm} \left\{
\abs{y(t, \xi, u)} - \gamma (\norm{u}) \right\}.
\ee
It follows from \rref{e-pf1} and \rref{e-pf10} that
\be{h-lower}
\abs{h(\xi)} \leq \bar{h} (\xi) \leq \beta_{0}
\left(\normx{\xi}\right),
\ee
where $\beta_{0} (s) = \beta (s, 0)$.
If $u_{0} \in \Mm$ such that $\beta_{0} (\normo{\xi}) \leq \gamma
(\norm{u_{0}})$, then
\beqn
& &\abs{y(t, \xi, u_{0})} - \gamma (\norm{u_{0}}) \le
\max \left\{ \beta \left(
\normo{\xi}, t \right), \gamma \left(
\norm{u_{0}} \right) \right\} - \gamma \left( \norm{u_{0}} \right)
= 0,\quad\forall\,t\ge 0.
\eeqn
Therefore, the supremum in \rref{1st-or} is achieved
when $\norm{u} \leq \gamma^{-1} \circ\beta_{0}
  (\norm{\xi})$, that is,
\be{2nd-or}
\bar{h} (\xi ) = \sup_{\substack{\norm{u} \leq \gamma^{-1}(
\beta_{0} (\normo{\xi})) \\ t \geq 0 }}
\left\{ \abs{y(t, \xi, u)} - \gamma (\norm{u}) \right\}.
\ee
Let $\tau \geq 0$ and $u \in \Mm$ be given, and
denote $\eta = x_{\tau}^{\xi, u}$. Then
\beq
& &\bar{h} \left( x_{\tau}^{\xi, u} \right) 
 =  \sup_{t \geq 0, v \in \Mm} \left\{
\abs{y(t, \eta, v)} - \gamma (\norm{v}) \right\}\nonumber \\
& = & \sup_{t \geq 0,  v\in \Mm} \left\{
\abs{y(t+\tau, \xi, \tilde{u})}
- \gamma (\norm{\tilde{u}}_{[\tau, \infty)}) \right\}\nonumber \\
& \leq & \sup_{t \geq 0, v \in \Mm} \left\{
\beta \left( \normo{\xi}, t+ \tau \right) + \gamma \left(
\norm{\tilde{u}} \right) - \gamma\left(\norm{\tilde{u}}_{[\tau, \infty)}\right)
\right\},\label{e-hbar2}
\eeq
where $\tilde{u}:= u \#_\tau v$
is defined by
\[
u \#_\tau v (t) = \left\{
\begin{array}{ll}
u(t) & \mbox{if $t \in [0, \tau)$}, \\
v(t-\tau) & \mbox{if $t \in [\tau, \infty)$.}
\end{array}
\right.
\]
Since
\be{triangle}
\gamma \left( \norm{\tilde{u}}_{[0, \infty)} \right)
\leq \gamma \left( \norm{\tilde{u}}_{[0, \tau)} \right)
+ \gamma \left( \norm{\tilde{u}}_{[\tau, \infty)} \right),
\ee
it follows from \rref{e-hbar2} that
\be{or-ios}
\bar{h} \left( x_{\tau}^{\xi, u} \right) \le \beta \left(\normo{\xi},
  \tau \right) + \gamma \left( \norm{\tilde{u}}_{[0, \tau)}
\right) \le
\beta \left(\normo{\xi},
  \tau\right) + \gamma \left( \norm{u}_{[0, \tau)}\right). 
\ee
This shows that the system with the output map $\bar h$ is \ios.

\bl{l-pr1}
System \rref{e-syso} with the output map $\bar h$ satisfies an
output Lagrange property as in \rref{ol}.
\els

\bpr
Similar to previous steps we consider
$\eta = x_{\tau}^{\xi, u}$ for $\tau \geq 0$ and  $u \in \Mm$. Then
\beqn
& &\bar{h} \left( x_{\tau}^{\xi, u} \right)
=  \sup_{t \geq 0, v \in \Uu} \left\{ \abs{y(t, \eta, v)}
- \gamma (\norm{v}) \right\} \\
& = &\sup_{t \geq 0, v \in \Mm} \left\{
 \abs{y(t+\tau, \xi, u\#_\tau v)}
 - \gamma (\norm{u\#_\tau v}_{[\tau, \infty)}) \right\} \\
&\leq&\!\!\!\sup_{s \geq 0, v \in \Mm} \left\{
\abs{y(s, \xi, \tilde{u})} - \gamma (\norm{\tilde u}_{[0, \infty)})
+ \gamma (\norm{\tilde u}_{[0, \tau)}) \right\} \\
&=& \bar{h} (\xi) + \gamma (\norm{u}) 
\le\max\{2\bar h(\xi), \ 2\gamma(\norm{u})\},
\eeqn
where $\tilde{u}:= u \#_\tau v$. This shows that an output Lagrange property \rref{ol} holds for
the output map $\bar h$.
\epr

\bl{l-l1} The map $\bar h$ is lower semi-continous.
\els

\bpr Let $\xi_0\in\xc$ and $\ve >0$. Denote $c= \bar{h} \left(
  \xi_{0} \right)$.
Then there exist $t_{0}$ and $u_{0}$ such that
\[
\abs{y(t_{0}, \xi_{0}, u_{0})}
- \gamma \left( \norm{u_0} \right) \geq c - \ve/2.
\]
By the continuity of the map $x_{t_0}^{\xi, u_0}$ in $\xi$ and the
continuity of $h$, there exists a neighborhood
$\mathcal B_{0} $ of $\xi_{0}$ such that
$$
\abs{ y(t_{0}, \xi, u_{0}) - y(t_{0}, \xi_{0}, u_{0}) }
\leq \frac{\ve}{2} \quad \forall \xi \in \mathcal B_{0}.
$$
Consequently,
\[
\bar{h} (\xi) \geq \abs{y(t_{0}, \xi, u_{0})}
- \gamma \left( \norm{u_{0}} \right) \geq c - \ve
\]
for every $\xi \in\mathcal B_{0}$. Hence,
\[
\liminf_{\xi \rightarrow \xi_{0}} \bar{h} (\xi) \geq \bar{h}
(\xi_{0}).
\]
This establishes the lower semi-continuity of $\bar{h}$.
\epr

To show the continuity (and the Lipschitz condition of $\bar h$ when $h$ is locally Lipschitz), we need the following result given under Lemma A.1 of
\cite{sw-ios-scl}: 
\bl{l-A1}
Let $\beta\in\kl$. Then there exists a family of mappings
$\{T_r\}_{r\ge 0}$ satisfying the following:
\begin{itemize}
\item for each fixed $r>0$, the map $T_r: \R_{>0}\rightarrow\R_{>0}$
  is continuous, onto, strictly decreasing, and $T_0(s)\equiv 0$;
\item for each fixed $s>0$, $T_r(s)$ is strictly increasing as $r$
  increases, and $\beta(r, t)<s$ for all $t\ge T_r(s)$.\qed
\end{itemize}
\els

Recall that $\Zz_{\bar h}=\{\xi\in\xc: \ \bar h(\xi) = 0\}$.
Pick
$\xi\notin \Zz_{\bar h}$.  Let $T_r(s)$ be as in Lemma \ref{l-A1} associated with $\beta$ given in \rref{e-pf1},
and let
$t_{\xi} = T_{\normx{\xi}} \left( \frac{\bar{h} (\xi)}{2}\right)$.
Then, 
\[
\abs{y(t, \xi, u)} -\gamma(\norm{u})
\le \beta(\normx{\xi}, t)\le \frac{\bar{h} (\xi)}{2}\qquad\forall\,t\ge t_\xi.
\]
It follows that
\be{fin-int-h-1}
\bar{h} (\xi ) = \sup_{\substack{\norm{u} \leq
\gamma^{-1} \left( \beta_{0} (\normo{\xi}) \right) \\
t \in [0, t_{\xi}] }}  \left\{
\abs{y(t, \xi, u)} - \gamma (\norm{u}) \right\}.
\ee
The significance of \rref{fin-int-h-1} is that, for $\xi\notin\Zz_{\bar h}$,
the supremum in the definition
of $\bar h$ is achieved in finite time proportional to $\bar h(\xi)$.

\bl{l-hbar-cont} The map $\bar h$ is continuous on $\xc\setminus\Zz_{\bar h}$.
\els

\bpr We first show that $\bar h$ is continuous on all compact subsets of $\xc\setminus\Zz_{\bar h}$.

Let $S\subseteq \xc\setminus\Zz_{\bar h}$ be compact, and let
$s_0>0$ such that  $\normx{\xi}\le s_0$ on $S$.
By the lower semi-continuity
of $\bar h$ and the compactness of $S$, there exists some
$c_0>0$ such that $\bar h(\xi)\ge c_0$
for all $\xi\in S$. Hence,
\be{fin-int-h-2}
\bar{h} (\xi) =
\sup_{\substack{\norm{u} \leq b \\ t \in [0, T] }}
\left\{ \abs{y(t, \xi, u)} - \gamma (\norm{u}) \right\}
\qquad \forall \xi \in S,
\ee
where $T = T_{s_0} \left(\frac{c_0}{2}\right)$ and
$b=\gamma^{-1} \left( \beta_{0} (s_{0}) \right)$.
Let
$$
{\mathscr C}_{S} = \left\{ \xu : t \in [0, T], \; \norm{u} \leq b, \;
\xi \in S \right\}.
$$

By Corollary \ref{c-lip}, the map $\xi\mapsto \xi_t^{\xi, u}$ is continuous on $S$, uniformly for $0\le t\le T$ and $u$ with $\norm{u}\le b$. Moreover, the set $\overline{{\mathscr C}_{S}}$ is compact, and hence, $h$ is uniformly continuous on the compact set $\overline{\mathscr C}_S$. As a consequence, for each $\ve > 0$, there exists $\delta>0$ such that
\be{e-conth}
\abs{h(x_t^{\xi, u})-h(x_t^{\eta, u})}< \ve,
\ee
for all $0\le t\le T$, all $\norm{u}\le b$, and all $\xi, \eta\in S$  such that $\normo{\xi-\eta}<\delta$.

Let $\ve>0$ be given.  Let $\delta>0$ be given as above.
For each $\xi\in S$, there exist
$t_{\xi, \ve} \in \left[0, T\right] $ and some
$u_{\xi, \ve}$ with $\norm{u_{\xi,\ve}}\le b$
such that
$$
\bar{h} (\xi )\leq \abs{y \left(
t_{\xi, \ve}, \xi, u_{\xi, \ve}\right)}
- \gamma \left( \norm{u_{\xi, \ve}} \right) + \ve.
$$
Consequently,
\beqn
\bar{h} (\xi) -\bar{h} (\eta)
&\leq& \abs{y \left(t_{\xi, \ve}, \xi,
u_{\xi, \ve}\right)} - \gamma \left(
\norm{u_{\xi, \ve}} \right)
+ \ve - \Big(\abs{y \left(t_{\xi, \ve}, \eta ,
u_{\xi, \ve}\right)}
- \gamma \left( \norm{u_{\xi, \ve } }\right) \Big) \\
&\leq&
\abs{h\left(x_{t_{\xi,\ve}}^{\xi, u_{\xi,\ve}}\right) - h\left(x_{t_{\xi,\ve}}^{\eta, u_{\xi,\ve}}\right)} +\ve < \ve + \ve
=2\ve,
\eeqn
for all $\xi, \eta\in S$ with $\normo{\xi-\eta} < \delta$.  By symmetry, it also holds that
$\bar{h} (\eta) -\bar{h} (\xi) < 2\ve$
if $\xi, \eta\in S$ and $\normo{\xi-\eta} < \delta$, and thus, $\abs{\bar h(\xi)-\bar h(\eta)} < 2\ve$ for all $\xi,\eta\in S$ with $\normo{\xi-\eta} < \delta$. This shows that $\bar h$ is continuous on every compact subset of $\xc\setminus\Zz_{\bar h}$.

To show the continuity of $\bar h$ on $\xc\setminus\Zz_{\bar h}$, it is enough to show the sequential continuity of $\bar h$.  Let $\xi_0\in\xc\setminus\Zz_{\bar h}$, and $\{\xi_k\}$ be a sequence converging to $\xi_0$.  Then the set $\{\xi_k: \ k\ge 1\}\cup\{\xi_0\}$ is compact.  It then follows that $\bar h(\xi_k)\rightarrow\bar h(\xi_0)$.  This proves the continuity of $\bar h$ on $\xc\setminus \Zz_h$.
\epr

\bl{l-l2} Assume that $h$ is locally Lipschitz, then the map $\bar h$ is locally Lipschitz on
$\xc\setminus \Zz_{\bar h}$.
\els

\bpr
We will show the local Lipschitz condition of $\bar h$ by showing that $\bar h$ is Lipschitz on compact subsets of
$\xc\setminus\Zz_{\bar h}$ (see Lemma \ref{l-Xu-Lip}.)  The proof will be a continuation of the proof of Lemma \ref{l-hbar-cont}. Still
let $S\subseteq \xc\setminus \Zz_{\bar h}$ be compact, and let
$s_0>0, c_0> 0, T>0,  b>0$, and the set ${\mathscr C}_{S}$ be
as in the proof of Lemma \ref{l-hbar-cont}.
By Lemma \ref{l-comp}, the set
$\overline{{\mathscr C}_{S}}$ is compact, and
thus, $h$ is Lipschitz on ${{\mathscr C}_{S}}$.  Furthermore,
the map $x_t^{\xi, u}$ is Lipschitz in $\xi\in S$, uniformly for $0\le t\le T$ and $\norm{u}\le b$.
It follows that
\[
\abs{y(t, \xi, u) - y(t, \eta, u)}\le
L\normo{x_t^{\xi, u} - x_t^{\eta, u}}\le LM\normo{\xi-\eta},
\]
for all $\xi,\eta\in\S$, all $t\in[0, T]$, and all $\norm{u}\le b$,
where
$L$ is Lipschitz constant of $h$ on ${\mathscr C}_S$, and $M$ is
the Lipschitz constant of $x_t^{\xi, u}$ in $\xi$ for $0\le t\le T$ and $\norm{u}\le b$.

For $\ve>0$, still let $t_{\xi,\ve}$ and $u_{\xi, \ve}$ be as in the proof of Lemma \ref{l-hbar-cont}.
Then,
\beqn
\bar{h} (\xi) -\bar{h} (\eta)
&\leq&  \abs{y \left(t_{\xi, \ve}, \xi,
u_{\xi, \ve}\right)} - \gamma \left(
\norm{u_{\xi, \ve}} \right)
+ \ve\\
& &\quad
- \left( \abs{y \left(t_{\xi, \ve}, \eta ,
u_{\xi, \ve}\right)}
- \gamma \left( \norm{u_{\xi, \ve } }\right) \right) \\
&\leq&
\abs{y \left(t_{\xi, \ve}, \xi, u_{\xi, \ve}\right)
- y \left(t_{\xi, \ve}, \eta , u_{\xi, \ve}\right)}
+ \ve \\
&\leq&
\hat L \bn \xi-\eta \enx  + \ve,
\eeqn
where $\hat L = LM$.
Note that the above is true for all $\ve>0$, it follows that
\[
\bar{h} (\xi) -\bar{h} (\eta) \le \hat L \bn \xi-\eta \enx,\quad
\forall\,\xi, \eta\in S.
\]
By symmetry of $\xi$ and $\eta$, one concludes that
\[
\abs{\bar{h} (\xi) -\bar{h} (\eta)} \le \hat L \bn \xi-\eta \enx,
\forall\, \xi, \eta\in S.
\]
This shows that $\bar h$ is Lipschitz on $S$.
\epr

\bl{l-l3} The map $\bar h$ is continuous on $\Zz_{\bar h}$.
\els

\bpr
Let $\xi_0\in \Zz_{\bar h}$. It is enough to show that for every sequence
$\{\xi_k\}\rightarrow \xi_0$, it holds that
$\{\bar h(\xi_k)\}\rightarrow 0$.

Suppose that this is not the case.  Then there exist a sequence
$\{\xi_k\}\rightarrow \xi_0$ and some $\ve_0>0$ such that $\bar
h(\xi_k)\ge \ve_0$ for all $k\ge 1$.  Without loss of generality, one
may assume that $\normx{\xi_k}\le s_1$ for all $k$, where $s_1>0$.  It
then follows that
\beqn
\bar h(\xi_k) = \sup_{0\le t\le t_2,\
\norm{u}\le b_2}\{\abs{y(t, \xi_k, u)} - \gamma(\norm{u})\},
\eeqn
where $t_2=T_{s_1}(\ve_0/2)$, and $b_2=\gamma^{-1}(\beta_0(s_1))$. It
then follows that for each $k$, there exist some $u_k$ with
$\norm{u_k}\le b_2$ and some $\tau_k\le t_2$ such that
\[
\abs{y(\tau_k, \xi_k, u_k)}-\gamma(\norm{u_k})
\ge \bar h(\xi_k) - \ve_0/2 \ge \ve_0/2.
\]
{Applying Corollary \ref{c-lip} to the compact set
$S_1:=\{\xi_k: \ k\ge 1\}\cup\{\xi_0\}$,} and by the continuity of
$h$, one sees that
\[
y(t, \xi_k, u)\rightarrow y(t, \xi_0, u) \quad
\mbox{as $k\rightarrow\infty$},
\]
uniformly in $t\in [0, t_2]$ and in $u$ with $\norm{u}\le b_2$.
This implies that
\beqn
\abs{y(\tau_k, \xi_0, u_k)} - \gamma(\norm{u_k})
\ge \abs{y(\tau_k, \xi_k, u_k)} - \ve_0/4 -
\gamma(\norm{u_k}) \ge \ve_0/4
\eeqn
for $k$ large enough, which contradicts the fact that $\bar
h(\xi_0)=0$. This proves the continuity of $\bar h$ on $\Zz_{\bar h}$.
Hence, one concludes the proof of Theorem \ref{olios-under-or}.
\epr

\subsection{Proof of \protect{Proposition \ref{p-robust}}}\label{s-62}

We will first prove Statement 2 followed by the proof of Statement 3, since the proof of Statement 1
depends on Statement 2 and Theorem \ref{olios-under-or}.

\subsubsection{Proofs of Statements 2 and 3}

As a preparation for a proof of Statement 1, we will prove the following version of Statement 2 of Proposition
\ref{p-robust} where $h$ is only assumed to be locally Lipschitz on $\xc\setminus \Zz_{h}$:

\bl{l-olios-olros-s}
Consider system \rref{e-syso} where $h:\xc\rightarrow \R^p$ is continuous everywhere and locally Lipschitz on $\xc\setminus \Zz_h$. Assume that the system is \rfc\ and is \olios.  Then there exists a smooth $\lambda\in\ki$ such that the corresponding system \rref{alt-sys} is \rfc\ and \olrgaos.
\el

Suppose system \rref{e-syso} is \rfc\ and \olios, and thus estimates \rref{e-rfc-e},
\rref{ios-max}
and \rref{ol} hold for some $\chi_1, \chi_2, \chi_3\in\kk$, $c\ge 0$,  $\beta\in\kl$, $\gamma\in\nn$, and
$\sigma\in\kk$. A function $d\in\MO$ will be referred as a disturbance for \rref{alt-sys}.
Without loss of generality, we assume that $\sigma(s)\ge s$,
and thus, $\sigma\in\ki$.  Furthermore, we assume that $\sigma(r)\ge \gamma(r)$ (if not, replace $\sigma(r)$ by $\max\{\sigma(r), \gamma(r)\}$.)


Let $\lambda \in\ki$ be smooth such that
\be{e-lambda}
\sigma(\lambda(\sigma(s)))< \frac s4, \quad\forall\,s>0
\ee
(which is equivalent to $\lambda(s) < \sigma^{-1}(\frac14\sigma^{-1}(s))$.) Since $\sigma(s)\ge s$, $\sigma(\lambda(s))\le s/4$.
To prove Proposition \ref{l-olios-olros-s}, we will establish several lemmas for system \rref{alt-sys} with the smooth map $\lambda\in\ki$ chosen as in \rref{e-lambda}.

First note that for the corresponding system \rref{alt-sys}, the map $g$ is locally Lipschitz on $\xc\setminus D$, where $D = \Zz_h$.
Hence, the uniqueness of trajectories cannot be assumed on $D$ in prior.

\bl{l-6.8} For any $\xi\notin D$ and $d\in\MO$, if $x_d(t, \xi, t)$ is defined on $[0, T)$ and $(x_d)_t^{\xi, d}\not\in D$ for all $t\in [0, T)$, then the following holds for all $t\in [ 0, T)$:
\be{e-33}
\sigma(\lambda(\abs{y_\lambda(t, \xi, d)}))\le\frac12\,\abs{h(\xi)}.
\ee
\els

\bpr
Let $\xi\in\xc$, $d\in\MO$, and $T>0$ be as given in the statement of the Lemma.  Then there is a unique solution of $\rref{alt-sys}$ starting from $\xi$ with the disturbance $d$ such that
$x_t^{\xi,d} \notin D$ on $[0, T)$.
Since
\[
\sigma(\lambda(\abs{h_\lambda(0, \xi, d)})) = \sigma(\lambda(\abs{h(\xi)}))\le\frac{\abs{h(\xi)}}{4},
\]
and since $h(\xi)\not=0$, it follows by continuity that
\[
\sigma(\lambda(\abs{h_\lambda(t, \xi, d)})) \le\frac12\abs{h(\xi)}
\]
for $t$ small enough.  Let
\[
t_1 = \inf\{t<T:\
\sigmay(\lambda(\abs{y_\lambda(t)}))>\frac12\abs{h(\xi)}\},
\]
(and $t_1=T$ if the set is empty.)  By continuity,
\be{e-t1}
\sigmay(\lambda(\abs{y_\lambda(t_1, \xi, d)}))=\frac12\abs{h(\xi)}.
\ee
Suppose $t_1 < T$.
Then \rref{e-33} holds on $[0, t_1)$. Consequently,
\be{e-l4-1}
\sigma\left(\abs{d(t)}\lambda(y_\lambda(t, \xi, d))\right)\le\frac12\abs{h(\xi)}
\ee
holds almost everywhere on $[0, t_1)$.
On $[0, t_1)$,
$y_\lambda(t)$ is an output function for \rref{e-syso} with the input
$u$ given by
$u(t) =d(t)\lambda(\abs{h(x_t)})$. By \rref{ol},
\beqn
\abs{y_\lambda(t)} &\le& \max\{\sigma(\abs{h(\xi)}), \ \sigma(
\abs{d(t)}\lambda(\abs{h((x_\lambda)_t)})\}
\le \sigma(\abs{h(\xi)})
\eeqn
on $[0, t_1)$, and thus,
\[
\sigma(\lambda(\abs{y_\lambda(t_1)})
\le \sigma(\lambda(\sigma(\abs{h(\xi)})))
\le \frac14\abs{h(\xi)}, 
\]
which contradicts \rref{e-t1}. This shows that
\rref{e-33} holds on $[0, T)$.
\epr

\bl{l-new-6}
The set $D$ is forward invariant for system \rref{alt-sys}, that is, for any $\xi\in D$ and any $d\in\MO$, the trajectory with the initial state $\xi$ and $d$ stays in $D$ on its maximum interval.~\mybox
\els

\bpr
Let $\xi\in D$ so that $h(\xi)=0$, and let $d\in\MO$.  Let $x_\lambda(t)$ be a trajectory with the initial state $\xi$ and the input $d$ defined on its maximum interval $[-\theta, t_{\max})$. For this case, we need to show that $y_\lambda(t)=0$ for all $t\in [0, t_{\max})$. Suppose this is not the case.  Then there exists some $0<t_1<t_{\max}$ and some $\ve>0$ such that
\[
\abs{y_\lambda(t_1)} \ge \ve.
\]
Let $\ve_0>0$ be such that $\ve_0 < \sigma(\lambda(\ve))$.  Then $\ve_0\le \ve/4$ (see \rref{e-lambda}).  Hence, there exists some $t_0\in (0, t_1)$ such that $\abs{y_\lambda(t_0)} = \ve_0$, and furthermore, one may assume that
\[
\abs{y_\lambda(t)}\ge \ve_0, \quad\forall t\in [t_0, t_1],
\]
(e.g., one may choose $t_0$ by $t_0:=\sup\{0\le t\le t_1: \ \abs{y_\lambda(t)}\le\ve_0\}$.)
Thus, on $[t_0, t_1]$, $(x_\lambda)_t\not\in D$.  Applying Lemma \ref{l-6.8} to the new initial state $\xi_1=(x_\lambda)_{t_0}$ 
one sees that
\beqn
\abs{y_\lambda(t)} \le \lambda^{-1}(\sigma^{-1}(
\abs{h(y_\lambda(t_0))}/2))\le \ve/2
\eeqn
for all $t\in [t_0, t_1]$. In particular, $\abs{y_\lambda(t_1)}
\le \ve/2$, which contradicts the assumption that
$\abs{y_\lambda(t_1)}\ge\ve$.  It follows that $y_\lambda(t)=0$
for all $t\in [0, t_{\max})$. This shows that the set $D$ is forward invariant.
\epr

\bc{c-uniq}(Uniqueness of Trajectories.) For each $\xi\in\xc$ and each $d\in\MO$, there is a unique trajectory of system \rref{alt-sys} corresponding to the initial state $\xi$ and the disturbance $d$ defined on its maximum interval.
\ec

\bpr
Let $\xi\in\xc$ and $d\in\MO$. If $\xi\in D$, then $(x_\lambda)_t\in D$ on its maximum interval, and thus, $d(t)\lambda(\abs{h((x_\lambda)_t)})\equiv 0$. Consequently, $x_\lambda(t)$ is a trajectory of the system
\be{e-fzero}
\dot x(t) = f_0(\xi):=f(x_t, 0).
\ee
Since $f_0$ is locally Lipschitz, there is a unique trajectory corresponding to the initial state $\xi$ and input $d$.

Assume that $\xi\not\in D$.  Let $\vf(t)$ and $\psi(t)$ be two trajectories of system \rref{alt-sys} corresponding to $d$, both satisfy the initial condition $\vf_0=\xi$ and $\psi_0=\xi$, defined on their maximum intervals $[-\theta, t_\vf)$ and $[-\theta, t_\psi)$ respectively. Since $g$ in \rref{alt-sys} is locally Lipschitz on $\xc\setminus D$, the uniqueness property holds on $\xc\setminus D$.

Suppose $\vf_{\bar t}\in D$ for some $\bar t\in [0, t_\vf)$.  Let
\[
\hat t = \inf\{t\le \bar t: \ \vf_t\in D\}.
\]
Then $\vf_{\hat t}\in D$ and $\vf_t\not\in D$ for all $0\le t<\hat t$. By the uniquness property on $\xc\setminus D$, it holds that
$\psi(t) = \vf(t)$ for all $t\in [0, \hat t)$.  By continuity, $h(\psi_{\hat t}) = h(\vf_{\hat t}) = 0$. By the uniqueness property on $D$, $\psi(t)=\vf(t)$ for all $t\ge \hat t$.  It then follows that $t_\vf=t_\psi$ and $\psi(t) = \vf(t)$ for all $t\in [-\theta, t_\vf)$.  Thus, for any $\xi\in\xc$ and any $d\in\MO$, there is a unique solution of \rref{alt-sys} corresponding to the initial state $\xi$ and the disturbance $d\in\MO$.
\epr

\bl{l-rfc-altsystem}
System \rref{alt-sys} is \rfc.
\el

\bpr
Let $\xi\in\xc$ and $d\in\MO$. By Lemmas \ref{l-6.8} and \ref{l-new-6},
\be{e-new33}
\sigma(\lambda(\abs{y_\lambda(t, \xi, d)}))\le \frac12 \abs{h(\xi)}
\ee
on the maximum interval $[-\theta, t_{\max})$ of the trajectory $x_\lambda(t, \xi, d)$.

By the uniqueness property of trajectories for \rref{e-syso} and \rref{alt-sys},
\be{e-twoeq}
x_\lambda(t, \xi, d) = x(t, \xi, u)\ {\rm with} \ u(t) = d(t)\lambda(|h(x_t^{\xi, u})|).
\ee
Consequently,
\[
\sigma(\abs{u(t)}) \le \sigma(\lambda(|h((x_\lambda)_t^{\xi, d})|)) \le  \frac12\abs{h(\xi)}, \quad {\rm a.e.}\ t\in [0, t_{\max}).
\]
Applying the \rfc\ property \rref{e-rfc-e} for system \rref{e-syso} with \rref{e-twoeq}, one sees that
\beq
\abs{x_\lambda(t, \xi, d)}
&\le& \chi_1(t) + \chi_2(\normo{\xi}) + \chi_3(\norm{u}) + c \nonumber\\
&\le & \chi_1(t) + \chi_2(\normo{\xi}) + \hat\chi_3(\abs{h(\xi)}/2) + c\nonumber\\
&\le & \chi_1(t) + \chi_4(\normo{\xi}) + c,\label{e-sysl-ufc}
\eeq
for all $t\in [0, t_{\max})$, where $\hat\chi_3(r) = \chi_3(\sigma^{-1}(r))$, and $\chi_4(r) = \chi_2(r) + \hat\chi_3(\frac{\pi(r)}2)$ (and where $\pi$ is as in \rref{bound-h}.)  It follows that $t_{\max}=\infty$, and \rref{e-sysl-ufc} holds on $[0, \infty)$.  This establishes the \rfc\ property for system \rref{alt-sys}.
\epr

A quick consequence of the \rfc\ property for \rref{alt-sys} is that \rref{e-new33} holds along all trajectories for all $t\ge 0$.
To complete the proof of Lemma \ref{l-olios-olros-s}, it remains to show that \rref{alt-sys} is \olrgaos.

Recall that the following \ios\ estimate holds for system \rref{e-syso} along all trajectories of \rref{e-syso}:
\be{e-iosmax1}
\abs{y(t, \xi, u)}\le \max\{\beta(\normo{\xi}, t), \ \sigma(\norm{u})\},\quad\forall\,t\ge 0.
\ee
Consequently, the following holds along all trajectories of \rref{alt-sys}:
\be{e-ios2}
\abs{y_\lambda(t, \xi, d)}\le \max\{\beta(\normo{\xi}, t), \ \sigma(\lambda(\abs{y_\lambda(t, \xi, d)}))\}
\le \max\left\{\beta(\normo{\xi}, t),\  \frac{\abs{h(\xi)}}{2}\right\}, 
\ee
for all $t\ge 0$.

\noindent{\it Claim.} For each $r>0$ and $s>0$, there exists some $T_{r, s}>0$
such that
\be{e-claim}
\abs{y_\lambda(t, \xi, d)}\le s/2,\qquad
\forall\,t\ge T_{r, s}, \ \forall\,\normo{\xi}\le r,\
\forall\, \abs{h(\xi)} \le s.
\ee

To prove the claim, we start with the \ios\ estimate \rref{e-ios2}:
\beqn
\abs{y_\lambda(t, \xi, d)} \le \max\left\{\beta(\normx{\xi}, t),\
\frac{\abs{h(\xi)}}{2}\right\}
\le\max\left\{\beta(r, t), \ \frac{s}{2}\right\}.
\eeqn
Since $\beta\in\kl$, there is some $T_{r,s} >0$ such that
$\beta(r, t) < s/2$ for all $t\ge T_{r, s}$.  It follows that
\[
\abs{y_\lambda(t, \xi, d)}\le \frac{s}{2}, \qquad
\forall\, t\ge T_{r, s}.
\]
The claim is thus proved.

Let $\ve >0$ be given. Pick any $\xi\not=0$ and $d$. Let
$r_1=\normo{\xi}$, and $s_1=\pi(r_1)$, where $\pi\in\ki$ is as in \rref{bound-h}. Note then that $\abs{h(\xi)}\le s_1$.  Let $s_{k+1}=\frac{s_{k}}2$ for $k\ge 1$.
Pick $L>0$ such that $2^{-L}s_1 < \ve$.

First of all, there is some $T_{r_1, s_1}>0$ such that
\[
\abs{y_\lambda(t, \xi, d)} \le s_1/2
\qquad\forall\, t\ge T_{r_1, s_1}.
\]
Let $T_1 = T_{r_1, s_1}$.
By the \rfc\ property \rref{e-sysl-ufc} (and assume that $\chi_1(r)\ge r$,) one has
\[
\normo{(x_\lambda)_{T_1}^{\xi,d}}
\le \chi_1(T_1) + \chi_4(r_1) + c.
\]
Let $r_2 = \chi_1(T_1) + \chi_4(r_1) + c$.
Applying \rref{e-claim} to $\normo{(x_\lambda)_{T_1}^{\xi,d}}\le r_2$ with
$\abs{y_\lambda(T_1, \xi, d)}\le s_1/2=s_2$, one sees that there
is some $T_{r_2, s_2}$ such that
\[
\abs{y_\lambda(t+T_{r_1, s_1}, \xi, d)}
=\abs{y_\lambda(t, \xi_1, d_1)}
\le \frac{s_2}{2} = \frac{s_1}{2^2},
\]
for all $t\ge T_{r_2, s_2}$, where $\xi_1=(x_\lambda)_{T_1}^{\xi, d}$, and
$d_1(t) = d(t+T_1)$.  Hence,
\[
\abs{y_\lambda(t, \xi, d)}\le \frac{s_2}{2}, \quad\forall\,t\ge T_2,
\]
where $T_2 = T_{r_1, s_1} + T_{r_2, s_2}$.

Inductively, assume that with $T_k = \sum_{i=1}^kT_{r_i, s_i}$, it holds that
\[
\abs{y_\lambda(T_k, \xi, d)}\le s_{k+1}.
\]
Applying
\rref{e-claim} to $y_\lambda(T_k, \xi, d)$ with
\[
r_{k+1} = \chi_1(T_k) + \chi_4(r_1) + c,
\]
one sees that there exists some $T_{r_{k+1}, s_{k+1}}$ so that
\[
\abs{y_\lambda(t+T_k, \xi, d)}\le \frac{s_{k+1}}{2}
= \frac{s_1}{2^{k+1}},
\]
for all $t\ge T_{r_{k+1}, s_{k+1}}$.  Let
\[
T = T_{r_1, s_1} + T_{r_2, s_2} + \cdots + T_{r_L, s_L}.
\]
Then, for all $t\ge T$,
\[
\abs{y_\lambda(t, \xi, d)} \le \frac{s_1}{2^L} < \ve.
\]
Note that $T$ depends only on $r_1$ and $\ve$.
Hence, it is shown that for any $\ve >0$
and $r>0$, there exists some $T_r(\ve)$ such that
\be{e-tre}
\abs{y_\lambda(t, \xi, d)} < \ve, \quad\forall\,t\ge T_r(\ve), \ \forall\,\normo{\xi}\le r.
\ee
Combining \rref{e-tre} with the output Lagrange property \rref{e-new33}, and noting that
$\abs{h(\xi)}\le\pi(\normo{\xi})$,  one concludes, by Lemma 7.5
in \cite{twj-dcdis} (applied with $a=\normo{\xi}$), that
there exists $\beta_1\in\kl$ such that the following holds along all trajectories of system \rref{alt-sys}:
\[
\abs{y_\lambda(t, \xi, d)}\le \beta_1(\normx{\xi}, t),\quad \forall\,t\ge 0.
\]
The proof of Lemma  \ref{l-olios-olros-s} is completed, so Statement 2 of Proposition \ref{p-robust} follows.

Below we continue to prove Statement 3, so we assume that \rref{e-rfc} and \rref{siios} hold for some $\chi_1, \chi_2, \chi_3\in\kk$, $c\ge 0$, $\beta\in\kk$ and $\gamma\in\Nn$.  
Without loss of generality, we rewrite \rref{siios} as
\be{e-siios-max}
\abs{y(t, \xi, u)}\le\max\{\beta(\abs{h(\xi)}, \ t), \ \gamma(\norm{u})\}, \quad\forall\,t\ge 0,
\ee
(by replacing $\beta$ and $\gamma$ with $2\beta$ and $2\gamma$ respectively.)
The idea of the proof will be the same as in the case of \olios.  Note that the \siios\ property \rref{siios} implies the \ol\ property: 
\[
\abs{y(t, \xi, u)} \le \max\{\sigma(\abs{h(\xi)}), \ \sigma(\norm{u}\},
\]
where $\sigma(r) = \max\{\beta(r, 0), \ \gamma(r)\}$.
Still define a smooth function $\lambda\in\ki$ as in \rref{e-lambda}.
Then Lemmas \ref{l-6.8}, \ref{l-new-6}, and \ref{l-rfc-altsystem} all still hold (and the uniqueness on $D$ is not an issue any more.)  It remains to show that the corresponding system \rref{alt-sys} is \sirgaos.  We follow the same steps as in the case of \olios, except now that property \rref{e-ios2} is strengthened to the following:
\be{e-siios-new2}
\abs{y_\lambda(t, \xi, d)}\le \max\left\{\beta(\abs{h(\xi)}, \ t), \ \frac{\abs{h(\xi)}}2\right\}, \quad\forall\,t\ge 0.
\ee
For each $r>0$, let $\Gamma_r>0$ be such that $\beta(r, t)\le \frac{r}{2}$ for all $t\ge \Gamma_r$. It then holds that if $\abs{h(\xi)}\le s$, then
\be{e-Ga}
\abs{y_\lambda(t, \xi, d)}\le \frac{s}{2}, \quad \forall\, t\ge \Gamma_r.
\ee
Let $\ve >0$ be given.  Pick any $\xi$ with $\abs{h(\xi)}\not=0$ and any $d\in\MO$, and let $r=\abs{h(\xi)}$.  Pick $L>0$ such that $2^{-L}r < \ve$. Let
\[
\hat\Gamma_r = \Gamma_r + \Gamma_{r/2} + \cdots + \Gamma_{r/(2^L)}.
\]
Repeatedly applying \rref{e-Ga}, one sees that for all $t\ge \hat\Gamma_r$
\[
\abs{y_\lambda(t, \xi, d)}\le \frac{r}{2^L} < \ve.
\]
Observe that $\hat\Gamma_r$ depends only on $r$ and $\ve$.  It follows that for any $\ve >0$, and any $r>0$, there exists some $T_r(\ve)$ such that
as long as $\abs{h(\xi)}\le r$, then
\[
\abs{y_\lambda(t, \xi, d)} < \ve, \qquad \forall\, t\ge T_r(\ve).
\]
Again, by Applying Lemma 7.5 in \cite{twj-dcdis} (this time with $a=\abs{h(\xi)}$), one concludes that there exists some $\hat\beta\in\kl$, such that
\[
\abs{y_\lambda(t, \xi, d)}\le \hat\beta(\abs{h(\xi)}, \ t),\quad\forall\,t\ge 0,
\]
along all trajectories of \rref{alt-sys}.  Statement 3 is thus proved.

\subsubsection{Proof of Statement 1}

The proof is based on Theorem
\ref{olios-under-or} and Lemma \ref{l-olios-olros-s}.


By Theorem \ref{olios-under-or}, there is an output-redefinition
$\bar h$, locally Lipschitz on $\xc\setminus \Zz_{\bar h}$, of $h$  such that the system \rref{e-syso} with $\bar h$ is \olios.
By Lemma \ref{l-olios-olros-s}, there exists some smooth $\lambda_0\in\ki$ such that the system
\be{e-ho}
\dot z(t) = f\left(z_t, d\lambda_0(\abs{\bar h(z_t)})\right),
\quad \bar y(t)=\bar h(z_t), 
\ee
is \rfc\ and \olrgaos, so for some $\beta_0\in\kl$, the following holds along all trajectories of \rref{e-ho}:
\be{e-ios-ho}
\abs{\bar y(t, \xi, d)}\le \beta_0(\normo{\xi}, t), \qquad \forall\,t\ge 0.
\ee
Recall that for some $\ab\in\ki$, it holds that $\ab(\abs{h(\xi)})\le |\bar h(\xi)|$.  In fact, from \rref{h-lower} one sees that $\ab(r) = r$ (so that one can always chooses $\ab(\cdot)$ to be smooth.) With $\lambda = \lambda_0\circ\ab$, it holds that
\[
\lambda(\abs{h(\xi)}) \le \lambda_0(\abs{\bar
h(\xi)}),\quad\forall\,\xi\in\xc.
\]

Below we show that the system
\be{e-12}
\dot x(t) = f(x_t, d\lambda(|h(x_t)|)),
\quad y(t)= h(x_t), 
\ee
where $d\in\MO$, is \rfc\ and satisfies the output convergence
property \rref{e-rgaos} for some $\beta\in\kl$.

Pick any $\xi$ and $d$. Let $\psi(\cdot)$ denote the
corresponding
trajectory of \rref{e-12}.
Suppose that $\psi(\cdot)$ is only defined on a maximum interval
$[0, T)$ for some $T<\infty$.
For $t\in [0, T)$, let
\[
d_0(t) = \frac{\lambda\left(\abs{h(\psi_t)}\right)}
{\lambda_0\left(\abs{\bar h(\psi_t)}\right)}\,d(t)\quad
\mbox{if $\bar h(\abs{\psi_t})\not=0$,}
\]
and $d_0(t)=0$ otherwise. Extend $d_0$ to
$[0, \infty)$ by letting $d_0(t) = 0$ for $t\ge T$. Then $d_0\in\MO$, and
\be{e-aa} 
\dot\psi(t) = f(\psi_t, d(t)\lambda(\abs{h(\psi_t)}))
= f(\psi_t, d_0(t)\lambda_0(\abs{\bar h(\psi_t)}))
\ee
on $[0, T)$.
Hence, $\psi(t)$ is a solution of the system \rref{e-ho} on $[0, T)$ with the
disturbance function $d_0$.
By the \rfc\ property of \rref{e-ho},
\be{e-ufc3}
\abs{\psi(t)}\le\chi_1(t) + \chi_2(\normo{\xi}) + c_1
\ee
on $[0, T)$ for some $\chi_1, \chi_2\in\kk$ and $c_1\ge 0$.
This contradicts to the
maximality of $T$.  Thus, $T=\infty$. In turn, this implies that
\rref{e-aa} holds for all $t\ge 0$, and
the \rfc\ estimate \rref{e-ufc3} holds on $[0, \infty)$ for
\rref{e-12}.  Furthermore,
\[
\abs{h(\psi_t)}\le \ab^{-1}\left(\abs{\bar
h(\psi_t)}\right)\le \beta(\normx{\xi}, t),
\]
where $\beta(s, t) =\ab^{-1}(\beta_0(s, t))$.  We have completed the proof that system \rref{e-12} is \rfc\ and \ios.
\qed

\medskip

\section{Conclusion}\label{s-7}

We extend the work in \cite{sw-ios-scl} on notions of input-to-output stability to delay systems by investigating
relations on several notions of output stability. For each of the \ios, \olios, and \siios\ properties, we provide trajectory-based Razumikhin criterion. The Lyapunov-based Razumikhin criteria for those properties will be developed in the forthcoming work. 
One of the main results, Theorem \ref{olios-under-or}, states that if a system is \ios, then the system is \olios\ under an output redefinition $\bar h$; and if the system map $f$ is locally Lipschitz, then $\bar h$ can be chosen to be locally Lipschitz on $\xc\setminus \Zz_{\bar h}$.  It still remains open if the local Lipschitz property of $\bar h$ can be extended to the whole space $\xc$.

Another question that remains open is that if the system map $f$ is further assumed to be Lipschitz on bounded subsets of $\xc$, then the output redefinition $\bar h$ in Theorem \ref{olios-under-or} can be chosen to be Lipschitz on bounded subsets of $\xc\setminus \Zz_{\bar h}$ or even on $\xc$.  This question is directly related to the question about the existence of Lyapunov-Krasovskii functional for the \ios\ property: if a system is \ios\ and if $f$ is assumed to be Lipschitz on bounded subsets of $\xc$, is there an associated Lyapunov-Krasovskii functional $V$ that is Lipschitz on bounded subsets of $\xc$?  For delay-free systems, the state space is a finite dimensional Euclidean space, where the local Lipschitz condition and the condition of being Lipschitz on bounded sets are equivalent. However, this equivalence fails on $\xc$. It is thus a natural question to understand if $\bar h$ can be chosen to be Lipschitz on bounded sets either on $\xc\setminus\Zz_{\bar h}$ or even on the whole space $\xc$ when the assumption on $f$ is strengthened to being Lipschitz on bounded sets.


\appendix

\section[*]{Appendix}

\subsection{Proofs of \protect{Lemma \ref{l-rfc}, Lemma \ref{l-comp},  and Corollary \ref{c-lip}}}\label{s-A1}

{\it Proof of Lemma \ref{l-rfc}.} The ``if'' implication trivially holds. The proof of the ``only if'' implication is almost the same as the proof of $(i)\Rightarrow(ii)$ in Lemma 2.12 in \cite{Mironchenko-Wirth-19}. To be more complete, we provide the proof to include the consideration of $u\in\Mm$.

Suppose that system \rref{e-syso} satisfies the \rfc\ property \rref{e-rfc}.
For each $r\ge 0$, let
\[
\rho(r) = \sup\{\abs{x(t, \xi, u)}: \ 0\le t\le r, \ \normo{\xi}\le r, \ \norm{u}\le r\}.
\]
By the \rfc\ property, $\rho(r) < \infty$ for each $r\ge 0$. Note that  $\rho(r)$ is nondecreasing. Let
\[
\rho_1(r) = \int_{r}^{r+1}\rho(s)\,ds.
\]
Then $\rho_1$ is continuous, nondecreasing, and $\rho_1(r)\ge \rho(r)$ for all $r\ge 0$.  
Let $c=\rho(0)$, and $\chi(r) = \rho_1(r) - \rho(0)$.  Then $\chi(0)=0$, and $\chi(\cdot)$ is continuous and $\rho(r) \le c+\chi(r)$.

Let $r_1, r_2, r_3 \ge 0$ be given.  By considering the three cases that $r_i=\max\{r_1, r_2, r_3\}$ for each $i\in\{1, 2, 3\}$, one sees
\[
\sup\{\abs{x(t, \xi, u)}: \ 0\le t\le r_1, \ \normo{\xi}\le r_2, \ \norm{u}\le r_3\}
\le \max\{\rho(r_1), \rho(r_2), \rho(r_3)\}.
\]
From this it follows that
\[
\abs{x(t, \xi, u)}\le \chi(t) + \chi(\normo{\xi}) + \chi(\norm{u}) + c,
\]
for all $t\ge 0$, $\xi\in\xc$, $u\in\Mm$.
\qed

\medskip

\noindent{\it Proof of Lemma \ref{l-comp}.}
Suppose system \rref{e-sysox} is \rfc.
Let $S\subseteq \xc$ be compact, $\Uu\subseteq\R^m$ be bounded, and
$T >0$.
It can be seen that the \rfc\ condition implies that $\Rr_{\Uu}^T(S)$ is bounded.
Let
\[
M = \max\{\abs{f(\xi, u)}: \ \xi\in\Rr_{\Uu}^T(S), \ u\in\Uu\}.
\]
By the
complete continuity of $f$, $M$ is well defined.  It then follows that
\[
\abs{x(t_1, \xi, u) - x(t_2, \xi, u)}\le M\abs{t_1 -t_2}
\]
for all $t_1, t_2\in [0, T], \xi\in S$ and $u\in\MmU$.

Since the set $S$ is compact, $S$ is equicontinuous. This means that
for every $\ve >0$, there is $\delta_0>0$ such that for all $\xi\in S$
\[
\abs{\xi(s_1) - \xi(s_2)} < \frac{\ve}{2}\quad\forall\,s_1, s_2\in [-\theta, 0]
\ {\rm and}\ \abs{s_1-s_2} < \delta_0.
\]
It can then be shown that for all $\xi\in S$, all $u\in\MmU$,
\[
\abs{x(t_1, \xi, u)-x(t_2, \xi, u)} < \ve, \quad
\forall t_1, t_2\in [-\theta, T]\ {\rm with}\ \abs{t_1-t_2} < \delta,
\]
where $\delta=\min\{\delta_0, \frac{\ve}{2M}\}$.
It then follows that the following holds for all $\xi\in S$, all $u\in\MmU$, and all $t\in[0, T]$:
\[
\normo{x_t^{\xi, u}(s_1) - x_t^{\xi, u}(s_2)} < \ve \quad \forall s_1, s_2\in [-\theta, 0]\ {\rm with}\ \abs{s_1-s_2} < \delta.
\]
Hence, the set
$\Rr_{\Uu}^T(S)$ is equicontinuous, and consequently, $\overline{\Rr_{\Uu}^T(S)}$ is bounded and equicontinuous.
By Arzel\'{a}-Ascoli Theorem, the set $\overline{\Rr_{\Uu}^T(S)}$ is compact.
\qed

\medskip

\noindent{\it Proof of Corollary \ref{c-lip}.}
To prove Part (a), let
$L$ be a Lipschitz constant of $f$ on the compact set $\overline{\Rr_{\Uu}^T(S)}$.  Then,
for all $\xi, \eta\in S$, $u\in\MmU$, and $0\le t\le T$:
\beq
\abs{x(t, \xi, u) - x(t, \eta, u)}\!
&\le& \abs{\xi(0)-\eta(0)} \nonumber\\ 
& &\quad {} +
\int_0^t\abs{f(x_\tau^{\xi, u}, u(\tau)) - f(x_\tau^{\eta, u}, u(\tau))}\,d\tau\nonumber\\
&\le&\!\abs{\xi(0)-\eta(0)} + L\int_0^t\normo{x_\tau^{\xi,u} - x_\tau^{\eta, u}}\,d\tau.
\label{e-flip3}
\eeq
For $t\in [-\theta, 0]$,
\[
\abs{x(t, \xi, u) - x(t, \eta, u)}\le\normo{\xi-\eta}.
\]
Thus, the following holds for all $t\in[-\theta, T]$, for all $\xi,\eta\in S$, and all $u\in\MmU$: 
\be{e-48}
\abs{x(t, \xi, u) - x(t, \eta, u)}
\le \normo{\xi-\eta} + L\int_0^t\normo{x_\tau^{\xi,u} - x_\tau^{\eta, u}}\,d\tau,
\ee
This implies that
for all $\xi, \eta\in S$, all $u\in\MmU$, and $s\in[-\theta, 0]$,
\[
\abs{x(t+s, \xi, u) - x(t+s, \eta, u)}\le \normo{\xi-\eta} + L\int_0^t\normo{x_\tau^{\xi,u} - x_\tau^{\eta, u}}\,d\tau,
\]
for all $t\in [0, T]$.  Taking the supremum for $s\in[-\theta, 0]$, one gets
\[
\normo{x_t^{\xi, u} - x_t^{\eta, u}}\le \normo{\xi-\eta} + L\int_0^t\normo{x_\tau^{\xi,u} - x_\tau^{\eta, u}}\,d\tau,
\]
for all $t\ge 0$. By Gr\"{o}nwall's inequality,
\[
\normo{x_t^{\xi, u} - x_t^{\eta, u}}\le \normo{\xi-\eta} e^{LT},
\]
for all $\xi, \eta\in S$, all $u\in\MmU$, and all $0\le t\le T$.

The above proof for Part (a) can easily be adapted to prove Part (b). Let $S$ and $\Uu$ be bounded, and $T>0$.
By the \rfc\ property, $\Rr_\Uu^T(S)$ is bounded.
Suppose $f$ is Lipschitz on bounded sets, then $f$ is Lipschitz on $\Rr_\Uu^T(S)$.  Hence one can still get \rref{e-flip3} and \rref{e-48}, and continue to
follow the same steps of Part (a) to complete the proof of Part (b).\qed

\subsection{Proof of \protect{Proposition \ref{2nd-form-olios}}}\label{s-olios}

In this section we include a proof for the Proposition \ref{2nd-form-olios} which provides a compact representation
of \olios\ that combines properties \rref{ios} and
\rref{ol} in a single estimate. This result was presented in \cite{sw-ios-scl}
and \cite{sw-ios-siam} but a detailed proof was skipped there. Though the complete proof can be found in \cite{Hasala-thesis},
to make the work more self-contained, we include the proof in this appendix.

The implication that an estimate as in \rref{e-olios2}
implies the \olios\ property is readily seen. The proof of the other
implication follows directly from the following lemma.
\bl{lagrange-conversion}
For any $\sigma \in \kk$ and $ \beta \in \kl$, there exist $ \hat{\rho} \in \kk$ and $ \hat{\beta} \in \kl$ such that,
\[
\min \left\{ \sigma (s), \beta (r,t)  \right\}
\leq \hat{\beta} \left( s, \frac{t}{1+ \hat\rho (r)}\right)
\]
for all $t, r, s \geq 0$.
\els

\bpr Let $\sigma\in\kk$ and $\beta\in\kl$.
By Proposition 7 in \cite{eds-iISS-scl98},
there exist $\at\in\ki$ and $\rho\in\ki$ such that, 
$\at \left(  \beta (r,t) \right) \leq  \rho (r) e^{-3t}$. For $r, s, t\ge 0$,
let
\[
a(r,s,t) := \min\left\{ \at \circ \sigma (s), \ \at \circ \beta (r,t)\right\}.
\]
Then $a(r,s,t) \le \sqrt{\at \circ \sigma (s)}\cdot \sqrt{\rho (r) e^{-3t}}$ for all $r, s, t\ge 0$.

Let $t, r, s \geq 0$ be given. First consider the case when
$t \geq \ln \left( 1 + \rho (r) \right)$.  In this case $\rho (r) e^{-t} \leq 1 $. Then,
\beqn
a(r,s,t) & \leq &  \sqrt{\at \circ \sigma (s)} \sqrt{\rho (r) e^{-3t}}
=
\sqrt{\at \circ \sigma (s)}\cdot e^{-t} \cdot \sqrt{\rho (r) e^{-t}} \\
&\leq& \sqrt{\at \circ
\sigma (s)}\cdot e^{-t} \leq \left( \sqrt{\at \circ \sigma (s) } \right) e^{\left(  -\frac{t}{1+\ln (1+ \rho (r))} \right) }.
\eeqn
Now consider the case when $t \in [0, \ln \left( 1 + \rho (r)\right)]$.
Then $t \leq 1 + \ln (1+ \rho (r))$ and hence;
\beqn
a(r,s,t)  &\leq&  \at \circ \sigma (s) \\
&=& \at \circ \sigma (s) \cdot e^{\left( - \frac{t}{1+\ln (1+ \rho (r))} \right) }
e^{\left( \frac{t}{1+\ln (1+ \rho (r))} \right)}\\
&\leq &  e\cdot\at \circ \sigma (s)  e^{\left( -\frac{t}{1+\ln (1+ \rho (r))} \right)}.
\eeqn
Let $ \tilde{\sigma} (s) = \max\{\sqrt{\at \circ \sigma (s) }, \
e\cdot\at \circ \sigma (s) \} $ for $s\in [0, \infty)$. Then
\[
a(r,s,t) \leq \tilde{\sigma} (s) e^{\left( -\frac{t}{1+\ln (1+ \rho (r))} \right)}\qquad\forall\,r, s, t\ge 0.
\]
It follows that
\[
\min \{ \sigma (s), \beta (r,t)\} \leq \hat{\beta} \left( s, \frac{t}{1+\hat{\rho}(r)} \right), \qquad\forall\, r, s, t\ge 0,
\]
where
$\hat{\beta} (s, \tau)  = \left(\at \right)^{-1} \left( \tilde{\sigma}(s) e^{-\tau} \right)$ and
$\hat{\rho}(r) = \ln (1+ \rho(r))$.
\epr

To complete the proof of Proposition \ref{2nd-form-olios}. 
suppose that system \rref{e-syso} is \olios.
Redefine $\gamma$ in \rref{ios} if necessary such that $\gamma(s) \geq \sigma_2 (s)$ for all $s\ge 0$,
where $\sigma_2 $ is as in \rref{ol}.
The case when $\norm{u}=\infty$ is trivial. Assume now that $\norm{u}<\infty$.
By Lemma \ref{lagrange-conversion}, there exist $\hat{\beta}
\in \kl$ and $\hat{\alpha} \in \kk$ satisfying:
\beqn
\abs{y(t, \xi, u) } - \gamma \left( \norm{u} \right)
&\leq& \min \{ \sigma_1 \left( \abs{h \left( \xi  \right)} \right) , \beta \left( \normx{\xi}, t \right)\}\\
&\leq& \hat{\beta} \left(\abs{h(\xi)}, \ \frac{t}{1+ \hat{\alpha} \left( \normx{\xi} \right)} \right), \qquad\forall\, t\ge 0.
\eeqn
Property \rref{e-olios2} thus follows.


\bibliographystyle{siam}
\bibliography{delayios}

\end{document}